\documentclass{article}
\usepackage{amssymb,amsmath,bm,hyperref}
\usepackage{amssymb,amsmath,bm,geometry,graphics,graphicx,url,color,inputenc}
\def\no{\noindent}
\def\pmatrix{\left(\begin{array}}
\def\endpmatrix{\end{array}\right)}
\def\udots{\reflectbox{$\ddots$}}  

\def\RR{\mathbb{R}}

\def\B{{\cal B}}

\def\I{{\cal I}}

\def\P{{\cal P}}

\def\dd{\mathrm{d}}

\newtheorem{theo}{Theorem}
\newtheorem{lem}{Lemma}
\newtheorem{cor}{Corollary}
\newtheorem{rem}{Remark}
\newtheorem{defi}{Definition}
\def\proof{\noindent\underline{Proof}\quad}
\def\QED{\mbox{~$\Box{~}$}}

\def\bfb{{\bm{b}}}
\def\bfc{{\bm{c}}}

\def\bfe{{\bm{e}}}

\def\bfzero{{\bm{0}}}
\def\bfuno{{\bm{1}}}

\def\bfgamma{{\bm{\gamma}}}

\def\bfeta{{\bm{\eta}}}
\def\bfvphi{{\bm{\varphi}}}

\def\bfpsi{{\bm{\psi}}}

\begin{document}

\title{High-order energy-conserving Line Integral Methods for charged particle dynamics}

\author{Luigi Brugnano\footnote{Dipartimento di Matematica e Informatica ``U.\,Dini'', Universit\`a di Firenze, Italy. \quad\url{luigi.brugnano@unifi.it}}  \and  Juan I.\,Montijano\footnote{I.U.M.A. -- Departamento de Matem\'atica Aplicada, Universidad de Zaragoza, Pza. San Francisco s/n, 50009 Zaragoza, Spain. \quad \url{{monti,randez}@unizar.es}} \and Luis R\'andez$^\dag$}


\maketitle

\begin{abstract} In this paper we study arbitrarily high-order  energy-conserving methods for simulating the dynamics of a charged particle.
They are derived and studied within the framework of {\em Line Integral Methods (LIMs)}, previously used for defining {\em Hamiltonian Boundary Value Methods (HBVMs)}, a class of energy-conserving Runge-Kutta methods for Hamiltonian problems. A complete analysis of the new methods is provided, which is confirmed by a few numerical tests. 

\medskip
\no{\bf Keywords:~} Charged particle dynamics, Lorentz force system, plasma physics, Energy-conserving methods, Line Integral Methods, LIMs, Hamiltonian Boundary Value Methods, HBVMs.

\medskip
\no{\bf MSC:~} 65L05, 65P10.

\end{abstract}

\section{Introduction}\label{intro} 
We shall here be concerned with the dynamics of a charged particle, which is described by the following system of ODEs,
\begin{equation}\label{cpd}
\dot q = p, \qquad \dot p = p\times L(q)-\nabla U(q), \qquad q(0)=q_0,~p(0)=p_0\in\RR^3,
\end{equation}
where $-\nabla U(q)$ and $L(q)$ are the electric and magnetic fields, respectively.\footnote{The equations (\ref{cpd}) are sometimes referred to as {\em Lorentz force system} (see, e.g., \cite{G1999,LW2016}).} 
Hereafter, we shall assume both of them to be time independent and suitably smooth functions of $q$.
For this motion, the energy,
\begin{equation}\label{H}
H(q,p) = \frac{1}2p^\top p+U(q),
\end{equation}
turns out to be conserved, along the solution of (\ref{cpd}) since, by the chain rule,
\begin{equation}\label{Hcons}
\frac{\dd}{\dd t} H(q,p) = \nabla U(q)^\top \dot q +  p^\top \dot p = p^\top\left(p\times L(q)\right) =0.
\end{equation}
Problem (\ref{cpd}) is a relevant one in plasma physics, since many important phenomena in plasmas can be understood and analyzed in terms of the motion of a single-particle \cite{Be2008}. This is due to the fact that collisions occur infrequently in hot plasmas and, in fact, (\ref{cpd}) provides a collisionless model of the plasma. However, despite the simplicity of the model, generic numerical methods are not appropriate for its long-time simulation, since an (unphysical) numerical drift in the energy (\ref{H}) is experienced, giving rise to a complete wrong solution orbit (see, e.g., \cite{QZXLST2013}).

Though several approaches have been recently proposed for numerically solving (\ref{cpd}) (see, e.g., \cite{HSLQ2015,HSLQ2016,HZSLQ2017,T2016,U2018,U2019}), a widely used method is the Boris method \cite{B1970}, which is symmetric and second-order accurate. In fact, because of its ease of implementation, such a method has become a {\em de facto} standard for simulating (\ref{cpd}). The analysis of its good behaviour, in turn, has been the subject of many investigations (see, e.g., \cite{EBQ2015,HL2018,QZXLST2013}). The Boris method, however, is not energy-conserving  and, in fact, a numerical drift may still be observed \cite{HL2018}.\footnote{In particular, the talk given by Ernst Hairer at the 2019 RSME Congress, presenting the results in \cite{HL2018}, has inspired the investigations reported in this paper.} On the other hand, a second-order energy conserving method based on a line-integral approach has been  recently proposed in \cite{LW2016}. Nevertheless, to the best of our knowledge no arbitrarily high-order energy-conserving methods are known for the simulation of problem (\ref{cpd}), and this motivates the present paper. The new energy-conserving methods will be derived within the framework of {\em Line Integral Methods} \cite{IaPa2007,IaPa2008,IaTr2009}, already used to devise a number of energy-conserving methods for various conservative problems \cite{ABI2015,BCMR2012,BI2012,BIT2012_1,BIT2012_2,BGI2018,BGIW2018}, with the main instance provided by Hamiltonian Boundary Value Methods (HBVMs), for the numerical solution of Hamiltonian problems \cite{BIT2009,BIT2010,BIT2012,BIT2015}. We also refer to the monograph \cite{LIMbook2016} and to the review paper \cite{BI2018}, for an overview.

The basic idea Line Integral Methods rely on is that of rewriting the conservation property (\ref{Hcons}) at $t=h$ in integral form, and defining  approximate paths $u(t)\approx q(t)$, and $v(t)\approx p(t)$ such that:\footnote{This procedure defines, indeed, the very first step of application of a one-step method.}
\begin{eqnarray}\label{lim1}
u(0) &=&q(0)~\equiv~ q_0, \qquad v(0) ~=~p(0)~\equiv~ p_0, \\[2mm]
\label{lim2}
u(h) &=:&q_1 ~\approx~q(h), \qquad v(h) ~=:~p_1 ~\approx~p(h), \\
\label{lim3}
H(q_1,p_1)-H(q_0,p_0) &\equiv&h\int_0^1\left[ \nabla U(u(ch))^\top\dot u(ch)+v(ch)^\top\dot v(ch)\right]\dd c ~=~0.
\end{eqnarray}
For our purposes, it will be convenient to rewrite (\ref{cpd}) in the equivalent form
\begin{equation}\label{cpd1}
\dot q = p, \qquad \dot p = B(q) p-\nabla U(q), \qquad q(0)=q_0,~p(0)=p_0\in\RR^3,
\end{equation}
where, setting $L(q)=\pmatrix{ccc} \ell_1(q),&\ell_2(q),&\ell_3(q)\endpmatrix^\top$,
\begin{equation}\label{Bq}
B(q) = \pmatrix{rrr} 0 &-\ell_3(q) &\ell_2(q)\\ \ell_3(q) &0 &-\ell_1(q)\\ -\ell_2(q) & \ell_1(q) &0\endpmatrix = -B(q)^\top.
\end{equation} 
In so doing, the arguments that we shall use for solving (\ref{cpd1})-(\ref{Bq}) can be naturally extended for solving more general problems, in the form 
\begin{equation}\label{general}
\ddot q = B(q) \dot q-\nabla U(q), \qquad q(0)=q_0,~\dot q(0)=p_0\in\RR^m, \qquad B(q)^\top = -B(q),
\end{equation}
which possess the same invariant (\ref{H}), with $\dot q=p$.

With this premise, the structure of the paper is as follows: in Section~\ref{new} we describe the framework in which the methods will be derived; in Section~\ref{discr} we provide a fully discrete method; in Section~\ref{discrete} its actual implementation is studied; in Section~\ref{num} we present a few numerical tests confirming the theoretical findings; at last, in Section~\ref{fine} we give some concluding remarks.

\section{Derivation of the method}\label{new}
Following the approach in \cite{BIT2012}, let us now rewrite the problem (\ref{cpd1})-(\ref{Bq}) by expanding the right-hand sides along a suitable orthonormal basis, which we choose as the  orthonormal Legendre polynomial basis $\{P_j\}$ on the interval $[0,1]$,
\begin{equation}\label{orto}
\deg P_i = i, \qquad \int_0^1 P_i(x)P_j(x)\dd x = \delta_{ij}, \qquad \forall i,j=0,1,\dots,
\end{equation}
with $\delta_{ij}$ the Kronecker symbol. We then obtain, at first: 
\begin{equation}\label{cpd2}
\dot q(ch) = \sum_{j\ge 0} P_j(c) \gamma_j(p), \qquad \dot p(ch) = \sum_{j\ge0} P_j(c)\left[B(q(ch)) \gamma_j(p)-\eta_j(q)\right],\qquad c\in[0,1],
\end{equation}
with 
\begin{equation}\label{gfe}
\gamma_j(p) = \int_0^1 P_j(\tau)p(\tau h)\dd\tau, \qquad
\eta_j(q) = \int_0^1 P_j(\tau)\nabla U(q(\tau h))\dd\tau.
\end{equation}
We observe that, in (\ref{cpd2}), $B(q(ch))$ has not yet been expanded. Next, consider the expansions
$$ 
P_j(c)B(q(ch)) = \sum_{i\ge0} P_i(c) \rho_{ij}(q), \qquad j=0,1,\dots,
$$
with 
\begin{equation}\label{ro}
\rho_{ij}(q)\equiv \rho_{ji}(q) =\int_0^1 P_i(\tau)P_j(\tau)B(q(\tau h))\dd\tau, \qquad i,j=0,1,\dots.
\end{equation}
As a result, from (\ref{cpd2}) we eventually arrive at:
\begin{equation}\label{cpd3}
\dot q(ch) = \sum_{i\ge 0} P_i(c) \gamma_i(p), \qquad \dot p(ch) = \sum_{i\ge0} P_i(c)\left[-\eta_i(q)+\sum_{j\ge0}\rho_{ij}(q) \gamma_j(p)\right],\qquad c\in[0,1].
\end{equation}

The following properties hold true.\footnote{Lemma~\ref{Ohj} is a generalization of \cite[Lemma\,1]{BIT2012}.}

\begin{lem}\label{Ohj} Assume $g:[0,h]\rightarrow V$, with $V$ a vector space, admit a Taylor expansion at 0. Then, for all $j=0,1,\ldots:$
$$\int_0^1 P_j(c)c^ig(ch)\dd c=O(h^{j-i}), \qquad i=0,\dots,j.$$\end{lem}
\proof 
By the hypotheses on $g$, one has: $$c^ig(ch)=\sum_{r\ge0} \frac{g^{(r)}(0)}{r!} h^r c^{r+i}.$$ 
Consequently, for all $i=0,\dots,j$, by virtue of (\ref{orto}) it follows that:
$$
\int_0^1 P_j(c)c^ig(ch)\dd c= \sum_{r\ge0} \frac{g^{(r)}(0)}{r!} h^r \int_0^1P_j(c)c^{r+i}\dd c
=\sum_{r\ge j-i} \frac{g^{(r)}(0)}{r!} h^r \int_0^1P_j(c)c^{r+i}\dd c = O(h^{j-i}).\QED 
$$               

\begin{cor}\label{grhj} With reference to (\ref{gfe}) and (\ref{ro}), for any suitably regular path  $\sigma:[0,h]\rightarrow\RR^3$  one has:\footnote{\label{foot5} The path would become $\sigma:[0,h]\rightarrow\RR^m$, in the case of the more general problem (\ref{general}).}  
\begin{equation}\label{gammaj}
\gamma_j(\sigma), \eta_j(\sigma)=O(h^j), \qquad \rho_{ij}(\sigma)=O(h^{|i-j|}).\qquad \forall i,j=0,1,\dots.
\end{equation}
\end{cor}
\bigskip

\begin{lem}\label{roijlem}
With reference to (\ref{ro}), for any path  $\sigma:[0,h]\rightarrow\RR^3$  one has:\footnote{ Footnote~\ref{foot5} applies also here.} 
\begin{equation}\label{roij}
\rho_{ij}(\sigma) = -\rho_{ij}(\sigma)^\top, \qquad \forall i,j=0,1,\dots.
\end{equation}
\end{lem} 
\proof The statement follows from the definition (\ref{ro}) and the skew-symmetry of matrix $B$ in (\ref{Bq}).\,\QED\bigskip

Next, in order to obtain polynomial approximations $u\approx q$ and $v\approx p$ of degree $s$, we truncate the infinite series in (\ref{cpd3}) after $s$ terms, thus getting:
\begin{equation}\label{u1v1}
\dot u(ch) = \sum_{i=0}^{s-1} P_i(c) \gamma_i(v), \qquad \dot v(ch) = \sum_{i=0}^{s-1} P_i(c)\left[-\eta_i(u)+\sum_{j=0}^{s-1}\rho_{ij}(u) \gamma_j(v)\right],\qquad c\in[0,1],
\end{equation}
with $\gamma_j(v), \eta_j(u), \rho_{ij}(u)$ defined according to (\ref{gfe}) and (\ref{ro}), by formally replacing $q,p$ with $u,v$, respectively.
 After that, we need to satisfy the {\em Line Integral conditions}  (\ref{lim1})--(\ref{lim3}), i.e.,
\begin{description}
\item[- requirement (\ref{lim1}):] we impose the initial conditions $u(0)=q_0$, $v(0)=p_0$, thus obtaining, integrating both equations in (\ref{u1v1}),
\begin{eqnarray}\label{uv}
u(ch) &=& q_0+h\sum_{i=0}^{s-1} \int_0^cP_i(x)\dd x\, \gamma_i(v), \\ \nonumber
v(ch) &=& p_0+h\sum_{i=0}^{s-1} \int_0^cP_i(x)\dd x\left[-\eta_i(u)+\sum_{j=0}^{s-1}\rho_{ij}(u) \gamma_j(v)\right], \qquad c\in[0,1];
\end{eqnarray}

\item[- requirement (\ref{lim2}):] accordingly, we set
\begin{equation}\label{uvh}
q_1:=u(h)\equiv q_0+h\gamma_0(v), \qquad p_1 := v(h) \equiv p_0-h\left[ \eta_0(u)-\sum_{j=0}^{s-1}\rho_{0j}(u)\gamma_j(v)\right],
\end{equation}
where, by virtue of (\ref{orto}), we took into account that $\int_0^1P_i(c)\dd c=\delta_{i0}$\,;

\item[- requirement (\ref{lim3}):] at last, next theorem states the property of energy-conservation.
\end{description}

\begin{theo}\label{H1H0} With reference to (\ref{H}) and (\ref{u1v1})--(\ref{uvh}), one has $H(q_1,p_1)=H(q_0,p_0)$.
\end{theo} 
\proof In fact, one has:
\begin{eqnarray*}
\lefteqn{H(q_1,p_1)-H(q_0,p_0) ~=~H(u(h),v(h))-H(u(0),v(0))~=~\int_0^h \frac{\dd}{\dd t} H(u(t),v(t))\dd t}\\
&=&h\int_0^1\left[ \nabla U(u(ch))^\top\dot u(ch)+v(ch)^\top\dot v(ch)\right]\dd c\\
&=&h\int_0^1\left[ \nabla U(u(ch))^\top\sum_{i=0}^{s-1}P_i(c)\gamma_i(v)\,+v(ch)^\top\sum_{i=0}^{s-1}P_i(c)\left(-\eta_i(u)+\sum_{j=0}^{s-1}\rho_{ij}(u)\gamma_j(v)\right)\right]\dd c\\
&=&h\sum_{i=0}^{s-1} \underbrace{\left[\int_0^1 P_i(c)\nabla U(u(ch))\dd c\right]^\top}_{=\,\eta_i(u)^\top}\gamma_i(v) -
h\sum_{i=0}^{s-1} \underbrace{\left[\int_0^1 P_i(c)v(ch)\dd c\right]^\top}_{=\,\gamma_i(v)^\top}\eta_i(v)\\
&&+\,h\sum_{i=0}^{s-1}\left[\int_0^1 P_i(c)v(ch)\dd c\right]^\top\sum_{j=0}^{s-1}\rho_{ij}(u)\gamma_j(u)\\
&=&h\sum_{i=0}^{s-1}\left[ \eta_i(u)^\top\gamma_i(v)-\gamma_i(v)^\top\eta_i(u)\right] + h\sum_{i,j=0}^{s-1}\gamma_i(v)^\top \rho_{ij}(u)\gamma_j(v)\\
&=&h\sum_{i,j=0}^{s-1}\gamma_i(v)^\top \rho_{ij}(u)\gamma_j(v)~=~0,
\end{eqnarray*}
where the last equality follows from (\ref{roij}).\,\QED
\bigskip

We now study to what extent $q_1,p_1$ approximate $q(h),p(h)$, respectively. For this purpose, in order to simplify the notation, we rewrite (\ref{cpd1})-(\ref{Bq}) as
\begin{equation}\label{y1}\dot y= f(y), \qquad y =\pmatrix{c} q\\ p\endpmatrix,\end{equation}
and denote by $y(t,t_0,y_0)$ the solution of such an equation, satisfying the initial condition $y(t_0)=y_0$. We also recall the following known perturbation results:
\begin{equation}\label{t0y0}
\frac{\partial}{\partial y_0} y(t,t_0,y_0) = \Phi(t,t_0), \qquad  \frac{\partial}{\partial t_0} y(t,t_0,y_0) = -\Phi(t,t_0)f(y_0),
\end{equation}
where $\Phi(t,t_0)$ is the fundamental matrix solution of the associated variational problem,
$$\dot\Phi(t,t_0) = f'(y(t,t_0,y_0))\Phi(t,t_0), \qquad \Phi(t_0,t_0)=I.$$
The following result then holds true, whose proof is based on the arguments used in \cite[Theorem~1]{BIT2012}.

\begin{theo}\label{y1yh} With reference to (\ref{cpd}) and (\ref{u1v1})--(\ref{uvh}), one has
 $$q_1=q(h)+O(h^{2s+1}), \qquad p_1=p(h)+O(h^{2s+1}),$$i.e., the approximation procedure has order $2s$.\end{theo} 
\proof In fact, by using the notation (\ref{y1}), the perturbation results (\ref{t0y0}), setting
$$w = \pmatrix{c} u\\ v\endpmatrix,\quad w_1 := w(h) \equiv\pmatrix{c}q_1\\ p_1\endpmatrix, \quad w_0 := w(0)\equiv\pmatrix{c}q_0\\ p_0\endpmatrix,\quad \Phi(t,t_0) \equiv \left[ \Phi_1(t,t_0),\, \Phi_2(t,t_0)\right],$$ 
with $\Phi_i(t,t_0)\in\RR^{6\times 3}$, $i=1,2$,\footnote{ $\Phi_i(t,t_0)\in\RR^{2m\times m}$ in the case of the more general problem (\ref{general}).} and considering that from (\ref{cpd2})--(\ref{cpd3})  (see also (\ref{u1v1}) and (\ref{uv})) it follows that
\begin{eqnarray*}
&&v(ch) = \sum_{j=0}^s P_j(c) \gamma_j(v),\qquad \nabla U(u(ch)) = \sum_{i\ge0} P_i(c)\eta_i(u), \\
 &&B(u(ch))v(ch) = \sum_{i\ge0} P_i(c) \sum_{j=0}^s \rho_{ij}(u) \gamma_j(v), \qquad~ c\in[0,1],
 \end{eqnarray*}
one has:
 \begin{eqnarray*}
 \lefteqn{
 w_1-y(h) ~=~w(h)-y(h) ~=~ y(h,h,w(h))-y(h,0,w(0)) ~=~ \int_0^h \frac{\dd}{\dd t} y(h,t,w(t))\dd t}\\
 &=&\int_0^h \left[\left.\frac{\partial}{\partial t_0} y(h,t_0,w(t))\right|_{t_0=t} + \left.\frac{\partial}{\partial y_0} y(h,t,y_0)\right|_{y_0=w(t)}
 \dot w(t)\right] \dd t\\
 &=& \int_0^h \left[-\Phi(h,t)f(w(t))+\Phi(h,t)\dot w(t)\right]\dd t ~=~ -h\int_0^1 \Phi(h,ch)\left[ f(w(ch))-\dot w(ch)\right]\dd c\\
 &=&-\,h\int_0^1 \Phi_1(h,ch)\left[ v(ch) - \sum_{i=0}^{s-1} P_i(c)\gamma_i(v)\right]\dd c\\
 &&+\,h\int_0^1 \Phi_2(h,ch)\left[\nabla U(u(ch))-B(u(ch))v(ch)- \sum_{i=0}^{s-1} P_i(c)\left(\eta_i(u)-\sum_{j=0}^{s-1}\rho_{ij}(u) \gamma_j(v)\right)\right]\dd c\\
 &=&-\,h \int_0^1 P_s(c)\Phi_1(h,ch)\dd c\,\gamma_s(v) + h\sum_{i\ge s}\int_0^1 P_i(c)\Phi_2(h,ch)\dd c\,\eta_i(u)\\
 &&-\,h\int_0^1 \Phi_2(h,ch)\left[\sum_{i=0}^{s-1} P_i(c) \rho_{is}(u) \gamma_s(v) +
 \sum_{i\ge s} P_i(c)\sum_{j=0}^s\rho_{ij}(u) \gamma_j(v)\right]\dd c \\
&\equiv& -\,h\Psi_{1s}\gamma_s(v) + h\sum_{i\ge s}\Psi_{2i}\eta_i(u)-h\sum_{i=0}^{s-1} \Psi_{2i}\,\rho_{is}(u) \gamma_s(v) - h \sum_{i\ge s} \Psi_{2i}\sum_{j=0}^s\rho_{ij}(u) \gamma_j(v) ~=:~(*),
 \end{eqnarray*}
 where we have set $$\Psi_{\ell\/i}:=\int_0^1 P_i(c)\Phi_\ell(h,ch)\dd c, \qquad \ell=1,2,\quad i\ge0.$$ We observe that, by virtue of Lemma~\ref{Ohj} and Corollary~\ref{grhj}, 
 $$\Psi_{\ell\/i},\gamma_i(v),\eta_i(u)=O(h^i),\qquad \rho_{ij}(u) = O(h^{|i-j|}), \qquad i,j\ge0.$$ 
 Consequently, one has:
 \begin{eqnarray*}
h\Psi_{1s}\gamma_s(v)=O(h^{2s+1}), &\qquad& h\sum_{i\ge s}\Psi_{2i}\eta_i(u)=O(h^{2s+1}),\\
\quad h\sum_{i=0}^{s-1}\underbrace{\Psi_{2i}\rho_{is}(u)}_{=\,O(h^s)} \gamma_s(v) = O(h^{2s+1}), 
&&h \sum_{i\ge s} \sum_{j=0}^s\Psi_{2i}\underbrace{\rho_{ij}(u) \gamma_j(v)}_{=\,O(h^i)}=O(h^{2s+1}),
\end{eqnarray*}
and, as a result, one concludes that $(*)=O(h^{2s+1})$.\,\QED\bigskip

 \section{Discretization}\label{discr}
 
 The approximation procedure (\ref{u1v1})--(\ref{uvh}) described in the previous section is not yet a {\em ready to use} numerical method. In fact, in order for this to happen, the integrals $\gamma_i(v), \eta_i(u), \rho_{ij}(u)$, $i,j=0,\dots,s-1$, defined in (\ref{gfe})-(\ref{ro}) need to be conveniently computed or approximated. We observe that, since $v\in\Pi_s$, then $\gamma_i(v)$ can be exactly computed by using the interpolatory quadrature formula of order $2s$ based at the zeros of $P_s$. If we denote $(\hat c_\ell,\hat b_\ell)$ the nodes and weights of such a quadrature,\footnote{I.e., $P_s(\hat c_\ell)=0$, $\ell=1,\dots,s$.} one has then:
\begin{equation}\label{hg}
\hat\gamma_i(v) := \sum_{\ell=1}^s \hat b_\ell P_i(\hat c_\ell) v(\hat c_\ell h) \equiv \gamma_i(v),\qquad i=0,\dots,s-1.
\end{equation}
We shall use the same quadrature for approximating $\rho_{ij}(u)$:
\begin{equation}\label{hro}
\hat\rho_{ij}(u) := \sum_{\ell=1}^s \hat b_\ell P_i(\hat c_\ell) P_j(\hat c_\ell) B(u(\hat c_\ell h)) \equiv \rho_{ij}(u)-\hat\Delta_{ij}(h), \qquad i,j=0.\dots,s-1,
\end{equation} with $\hat\Delta_{ij}(h)$ the quadrature error. 
Finally, as it has been done in the case of HBVMs \cite{BIT2010},   for approximating $\eta_i(v)$ we shall use a Gauss-Legendre quadrature of order $2k$, i.e., the interpolatory quadrature rule based at the zeros of $P_k$, for a convenient value $k\ge s$, with abscissae and weights $(c_\ell,b_\ell)$:\footnote{I.e., $P_k(c_\ell)=0$, $\ell=1,\dots,k$.}
\begin{equation}\label{he}
\hat\eta_i(u) := \sum_{\ell=1}^k b_\ell P_i(c_\ell) \nabla U(c_\ell h) \equiv \eta_i(u)-\Delta_i(h), \qquad i=0,\dots,s-1,
\end{equation} with $\Delta_i(h)$ the quadrature error.  Concerning the quadrature errors, one verifies that, for all $i,j=0,\dots,s-1:$
\begin{equation}\label{Delta}
\hat\Delta_{ij}(h) = O(h^{2s-i-j}), \qquad \Delta_i(h) = \left\{ \begin{array}{cc} 0, &\mbox{if}~U\in\Pi_\nu~\mbox{with}~\nu\le 2k/s,\\[3mm] O(h^{2k-i}), &\mbox{otherwise.}\end{array}\right.
\end{equation} 
Consequently, it is straightforward to prove that the results of Corollary~\ref{grhj} and Lemma~\ref{roijlem} continue formally to hold for $\hat\gamma_i(v), \hat\eta_i(u), \hat\rho_{ij}(u)$ defined in (\ref{hg})--(\ref{he}) (and, of course, $\hat\rho_{ij}(u)=\hat\rho_{ji}(u)$).
In so doing, the polynomials (\ref{u1v1})--(\ref{uv}) respectively become\,\footnote{For sake of brevity, we shall continue to denote such polynomials by $u$ and $v$, respectively.}
\begin{eqnarray}\label{newuv1}
\dot u(ch) &=& \sum_{i=0}^{s-1} P_i(c) \hat\gamma_i(v), \\ \nonumber
\dot v(ch) &=& \sum_{i=0}^{s-1} P_i(c)\left[-\hat\eta_i(u)+\sum_{j=0}^{s-1}\hat\rho_{ij}(u) \hat\gamma_j(v)\right],\qquad c\in[0,1],
\end{eqnarray}
and
\begin{eqnarray} \label{newuv}
u(ch) &=& q_0+h\sum_{j=0}^{s-1} \int_0^cP_j(x)\dd x\, \hat\gamma_j(v), \\ \nonumber
v(ch) &=& p_0+h\sum_{i=0}^{s-1} \int_0^cP_i(x)\dd x\left[-\hat\eta_i(u)+\sum_{j=0}^{s-1}\hat\rho_{ij}(u) \hat\gamma_j(v)\right], \qquad c\in[0,1],
\end{eqnarray}
with the new approximations given by
\begin{equation}\label{newuvh}
q_1:=u(h)\equiv q_0+h\hat\gamma_0(v), \qquad
p_1 := v(h) \equiv p_0-h\left[ \hat\eta_0(u)-\sum_{j=0}^{s-1}\hat\rho_{0j}(u)\hat\gamma_j(v)\right],
\end{equation}
in place of (\ref{uvh}). Clearly, the {\em Line Integral requirements} (\ref{lim1}) and (\ref{lim2}) are satisfied by the new polynomial paths $u$ and $v$ defined by (\ref{hg})--(\ref{newuvh}). Concerning the requirement (\ref{lim3}), the following result holds true.

\begin{theo}\label{newH1H0} With reference to (\ref{H}) and (\ref{hg})--(\ref{newuvh}), one has either $H(q_1,p_1)=H(q_0,p_0)$, if\, $U\in\Pi_\nu$ with $\nu\le 2k/s$, or $H(q_1,p_1)=H(q_0,p_0)+O(h^{2k+1})$, otherwise.
\end{theo} 
\proof In fact, one has, by considering that \,$\gamma_i(v)=\hat\gamma_i(v)$\, and \,$\hat\rho_{ij}(u) = \hat\rho_{ji}(u) =-\hat\rho_{ij}(u)^\top$:
\begin{eqnarray*}
\lefteqn{H(q_1,p_1)-H(q_0,p_0) ~=~H(u(h),v(h))-H(u(0),v(0))~=~\int_0^h \frac{\dd}{\dd t} H(u(t),v(t))\dd t}\\
&=&h\int_0^1\left[ \nabla U(u(ch))^\top\dot u(ch)+v(ch)^\top\dot v(ch)\right]\dd c\\
&=&h\int_0^1\left[ \nabla U(u(ch))^\top\sum_{i=0}^{s-1}P_i(c)\gamma_i(v)\,+v(ch)^\top\sum_{i=0}^{s-1}P_i(c)\left(-\hat\eta_i(u)+\sum_{j=0}^{s-1}\hat\rho_{ij}(u)\gamma_j(v)\right)\right]\dd c\\
&=&h\sum_{i=0}^{s-1} \underbrace{\left[\int_0^1 P_i(c)\nabla U(u(ch))\dd c\right]^\top}_{=\,\eta_i(u)^\top}\gamma_i(v) -
h\sum_{i=0}^{s-1} \underbrace{\left[\int_0^1 P_i(c)v(ch)\dd c\right]^\top}_{=\,\gamma_i(v)^\top}\hat\eta_i(v)\\
&&+h\sum_{i=0}^{s-1}\left[\int_0^1 P_i(c)v(ch)\dd c\right]^\top\sum_{j=0}^{s-1}\hat\rho_{ij}(u)\gamma_j(u)\\
\end{eqnarray*}\begin{eqnarray*} 
&=&h\sum_{i=0}^{s-1}\left[ \eta_i(u)^\top\gamma_i(v)-\gamma_i(v)^\top\left(\eta_i(u)-\Delta_i(h)\right)\right] + h\sum_{i,j=0}^{s-1}\gamma_i(v)^\top \hat\rho_{ij}(u)\gamma_j(v)\\
&=&h\sum_{i=0}^{s-1}\gamma_i(v)^\top \Delta_i(u).
\end{eqnarray*}
The statement then follows from (\ref{Delta}) and considering that $\gamma_i(u)=O(h^i)$.\,\QED
\bigskip

\begin{rem}\label{klarge} As is clear from Theorem~\ref{newH1H0}, an exact energy conservation is obtained in the polynomial case, by choosing $k$ large enough. However, also in the non-polynomial case, one can always gain a {\em practical} energy conservation, by choosing $k$ large enough so that the $O(h^{2k+1})$ energy error falls within the round-off error level.\end{rem} 

Next results states that the order $2s$ of the approximation procedure (\ref{u1v1})--(\ref{uvh}) is retained by the new one.

\begin{theo}\label{newy1yh} With reference to (\ref{cpd}) and (\ref{hg})--(\ref{newuvh}), for all $k\ge s$ one has
 $$q_1=q(h)+O(h^{2s+1}), \qquad p_1=p(h)+O(h^{2s+1}).$$\end{theo} 
\proof In fact, by using the same notations and preliminary results used in the proof of Theorem~\ref{y1yh},
one has:
 \begin{eqnarray*}
 \lefteqn{
 w_1-y(h) ~=~w(h)-y(h) ~=~ y(h,h,w(h))-y(h,0,w(0)) ~=~ \int_0^h \frac{\dd}{\dd t} y(h,t,w(t))\dd t}\\
 &=&\int_0^h \left[\left.\frac{\partial}{\partial t_0} y(h,t_0,w(t))\right|_{t_0=t} + \left.\frac{\partial}{\partial y_0} y(h,t,y_0)\right|_{y_0=w(t)}
 \dot w(t)\right] \dd t\\
 &=& \int_0^h \left[-\Phi(h,t)f(w(t))+\Phi(h,t)\dot w(t)\right]\dd t ~=~ -h\int_0^1 \Phi(h,ch)\left[ f(w(ch))-\dot w(ch)\right]\dd c\\
 &=&-\,h\int_0^1 \Phi_1(h,ch)\left[ v(ch) - \sum_{i=0}^{s-1} P_i(c)\gamma_i(v)\right]\dd c\\
 &&+\,h\int_0^1 \Phi_2(h,ch)\left[\nabla U(u(ch))-B(u(ch))v(ch)- \sum_{i=0}^{s-1} P_i(c)\left(\hat\eta_i(u)-\sum_{j=0}^{s-1}\hat\rho_{ij}(u) \gamma_j(v)\right)\right]\dd c\\
 &=&-\,h\int_0^1 \Phi_1(h,ch)\left[ v(ch) - \sum_{i=0}^{s-1} P_i(c)\gamma_i(v)\right]\dd c +h\int_0^1 \Phi_2(h,ch)
 \left[\nabla U(u(ch))-B(u(ch))v(ch) \begin{array}{c} ~\\~\\~ \end{array} \right.\\ &&\left.
 - \sum_{i=0}^{s-1} P_i(c) \left(
 \left(\eta_i(u) -\Delta_i(h)\right)  -\sum_{j=0}^{s-1}\left(\rho_{ij}(u)-\hat\Delta_{ij}(h)\right) \gamma_j(v)\right)\right]\dd c\\
 &=& O(h^{2s+1}) + h\sum_{i=0}^{s-1} \Psi_{2i}\Delta_i(h) -h\sum_{i,j=0}^{s-1} \Psi_{2i}\hat\Delta_{ij}(h)\gamma_j(v),
\end{eqnarray*}
 where the last equality follows from the proof of Theorem~\ref{y1yh}, in which formally the same terms were involved, except those including the quadrature errors (\ref{Delta}). Concerning these latter terms,  by recalling that
  $$\Psi_{2i},\gamma_i(v)=O(h^i),\qquad \Delta_i(h)=O(h^{2k-i}),\qquad  \hat\Delta_{ij}(h) = O(h^{2s-i-j}), \qquad i,j=0,\dots,s-1,$$ 
 one has:
$$h\sum_{i=0}^{s-1} \Psi_{2i}\Delta_i(h) = O(h^{2k+1}), \qquad h\sum_{i,j=0}^{s-1} \Psi_{2i}\hat\Delta_{ij}(h) \gamma_j(v) = O(h^{2s+1}).$$
Consequently, the statement follows.\,\QED\bigskip

\begin{defi}\label{lim} Hereafter, we shall refer to the method defined by (\ref{hg})--(\ref{newuvh}) as {\em Line Integral Method with parameters $(k,s)$}, in short {\em LIM$(k,s)$}, for solving problem (\ref{cpd1})--(\ref{Bq}).\end{defi} 

We then conclude that, according to Theorem~\ref{newy1yh}, the LIM$(k,s)$ has order $2s$, for all  allowed values of $s$ and $k\ge s$. Moreover, according to Theorem~\ref{newH1H0} and Remark~\ref{klarge}, it is energy-conserving, either exactly or {\em practically}, by choosing $k$ suitably large. This, in turn, will be not a drawback, since the discrete problem that one has to solve  has dimension $s$, {\em independently of $k$}, as we are going to see in the next section.

\section{The discrete problem}\label{discrete}
We now study the efficient implementation of the method (\ref{hg})--(\ref{newuvh}). For this purpose, we need the following matrices, defined by the Legendre polynomial basis (\ref{orto}) and the nodes and weights of the Gauss-Legendre quadratures of order $2s$, $(\hat c_\ell,\hat b_\ell)$, and $2k$, $(c_\ell,b_\ell)$:
\begin{eqnarray}\nonumber
\hat\P_s =\pmatrix{ccc} P_0(\hat c_1)&\dots &P_{s-1}(\hat c_1)\\ \vdots & &\vdots\\ P_0(\hat c_s)&\dots &P_{s-1}(\hat c_s)\endpmatrix,&& 
\hat\I_s =\pmatrix{ccc} \int_0^{\hat c_1}P_0(x)\dd x&\dots &\int_0^{\hat c_1}P_{s-1}(x)\dd x\\ \vdots & &\vdots\\ 
\int_0^{\hat c_s}P_0(x)\dd x&\dots &\int_0^{\hat c_s}P_{s-1}(x)\dd x\endpmatrix,\\[3mm] \label{PIO}
\P_s =\pmatrix{ccc} P_0(c_1)&\dots &P_{s-1}(c_1)\\ \vdots & &\vdots\\ P_0(c_k)&\dots &P_{s-1}(c_k)\endpmatrix,&&
\I_s =\pmatrix{ccc} \int_0^{c_1}P_0(x)\dd x&\dots &\int_0^{c_1}P_{s-1}(x)\dd x\\ \vdots & &\vdots\\ 
\int_0^{c_k}P_0(x)\dd x&\dots &\int_0^{c_k}P_{s-1}(x)\dd x\endpmatrix,\\[3mm] \nonumber
\hat\Omega = \pmatrix{ccc} \hat b_1\\ &\ddots\\ &&\hat b_s\endpmatrix, &&
\Omega = \pmatrix{ccc} b_1\\ &\ddots\\ &&b_k\endpmatrix.
\end{eqnarray}   
Moreover, we need to introduce the vector \,$\bfuno_r=\pmatrix{ccc}1,&\dots,1\endpmatrix^\top\in\RR^r$, and the block vectors and matrices, with reference to (\ref{hg})--(\ref{he}):
\begin{equation}\label{geG}
\bfgamma(v) := \pmatrix{c} \hat\gamma_0(v)\\ \vdots \\ \hat\gamma_{s-1}(v)\endpmatrix,\,
\bfeta(u) := \pmatrix{c} \hat\eta_0(u)\\ \vdots \\ \hat\eta_{s-1}(u)\endpmatrix,\,
\Gamma(u) := \pmatrix{ccc} \hat\rho_{00}(u)& \dots &\hat\rho_{0,s-1}(u)\\ \vdots & &\vdots\\
\hat\rho_{s-1,0}(u)& \dots &\hat\rho_{s-1,s-1}(u)\endpmatrix.
\end{equation}

We observe that, for computing what in (\ref{geG}), we need to evaluate $v(\hat c_\ell h), u(\hat c_\ell h)$, $\ell = 1,\dots,s$, and $u(c_\ell h)$, $\ell=1,\dots,k$. For this purpose, to have a more compact notation, we define the vectors:
$$\hat\bfc := \pmatrix{c} \hat c_1\\ \vdots \\ \hat c_s\endpmatrix,
\bfc := \pmatrix{c} c_1\\ \vdots \\ c_k\endpmatrix,
v(\hat\bfc h):= \pmatrix{c} v(\hat c_1h)\\ \vdots \\ v(\hat c_s h)\endpmatrix,
u(\hat\bfc h):= \pmatrix{c} u(\hat c_1h)\\ \vdots \\ u(\hat c_s h)\endpmatrix,
u(\bfc h):= \pmatrix{c} u(c_1h)\\ \vdots \\ u(c_k h)\endpmatrix.$$
A similar notation will be used for the functions in (\ref{cpd})--(\ref{Bq}), when a block vector argument is specified.
With reference to (\ref{PIO}) and (\ref{geG}), one then obtains
\begin{equation}\label{uhc}
u(\hat\bfc h) ~= \bfuno_s\otimes q_0 + h\hat\I_s\otimes I\bfgamma(v),\qquad
u(\bfc h)      ~=~ \bfuno_k\otimes q_0 + h\I_s\otimes I\bfgamma(v),
\end{equation}
and
$$
v(\hat\bfc h) ~=~ \bfuno_s\otimes p_0 + h\hat\I_s\otimes I\left( \Gamma(u)\bfgamma(v) -\bfeta(u)\right).
$$ Hereafter, $I$ will denote the identity matrix having the size of the vectors $q_0$ and $p_0$ (i.e., 3 for problem (\ref{cpd}),  or $m$ for problem (\ref{general})).
The last equation can be rewritten, by defining the block vector
\begin{equation}\label{bfi}
\bfvphi(u,v)~\equiv~\pmatrix{c} \varphi_0(u,v)\\ \vdots \\ \varphi_{s-1}(u,v)\endpmatrix ~:=~ \Gamma(u)\bfgamma(v) ,
\end{equation} 
as:
\begin{equation}\label{vhc}
v(\hat\bfc h) ~=~ \bfuno_s\otimes p_0 + h\hat\I_s\otimes I\left( \bfvphi(u,v) -\bfeta(u)\right).
\end{equation}
From (\ref{newuvh}), one then obtains that the new approximations are given by
\begin{equation}\label{newuh_1}
q_1 = q_0 + h\hat\gamma_0(v), \qquad p_1 = p_0 + h\left[ \varphi_0(u,v) - \hat\eta_0(u)\right].
\end{equation}
Consequently, we need to compute the three block vectors $\hat\bfgamma(v)$, $\hat\bfeta(u)$, and $\bfvphi(u,v)$, respectively defined in (\ref{geG}) and (\ref{bfi}). By using the matrices defined in (\ref{PIO}), and the formulae (\ref{hg})--(\ref{he}), we then obtain
\begin{equation}\label{bfge}
\bfgamma(v) = \hat\P_s^\top \hat\Omega\otimes I\, v(\hat\bfc h), \qquad \bfeta(u) = \P_s^\top \Omega\otimes I\,\nabla U( u(\bfc h)),
\end{equation}
and
\begin{equation}\label{bffi1}
\bfvphi(u,v) = \hat\P_s^\top \hat\Omega\otimes I\, \B( u(\hat\bfc h) ) \hat\P_s\otimes I\,\bfgamma(v),
\end{equation}
having set
$$\B(u(\hat\bfc h)) = \pmatrix{ccc} B(u(\hat c_1h))\\ &\ddots\\ &&B(u(\hat c_s h))\endpmatrix.$$
Nevertheless, by taking into account that $\hat\P_s^\top\hat\Omega = \hat\P_s^{-1}$  and the first equation in (\ref{bfge}), one obtains that (\ref{bffi1}) can be rewritten as
\begin{equation}\label{bffi2}
\bfvphi(u,v) = \hat\P_s^\top \hat\Omega\otimes I\, \left[\B( u(\hat\bfc h) ) v(\hat\bfc h)\right].
\end{equation}
Finally, by considering the expressions (\ref{uhc}) and (\ref{vhc}), we eventually obtain the following set of equations:

\begin{eqnarray}\label{gam}
\bfgamma &=& \hat\P_s^\top \hat\Omega\otimes I\, \left[ \bfuno_s\otimes p_0 + h\hat\I_s\otimes I\left( \bfvphi -\bfeta\right)\right],\\
\label{eta}
\bfeta &=& \P_s^\top \Omega\otimes I\,\nabla U\left( \bfuno_k\otimes q_0 + h\I_s\otimes I\bfgamma \right),\\ \label{vfi}
\bfvphi &=& \hat\P_s^\top \hat\Omega\otimes I\, \left[\B\left( \bfuno_s\otimes q_0 + h\hat\I_s\otimes I\bfgamma \right)\cdot\left(\bfuno_s\otimes p_0 + h\hat\I_s\otimes I\left( \bfvphi -\bfeta\right)\right)\right],
\end{eqnarray}
where we have removed the arguments of the (block) vectors, since now this is only an algebraic system of equations which, moreover,
 can be further simplified.  In fact, we observe, at first, by taking into account
\begin{equation}\label{e1}
\hat\P_s^\top\Omega \bfuno_s = \pmatrix{cccc} 1,& 0,&\dots,&0\endpmatrix^\top =:\bfe_1 \in\RR^s, 
\end{equation}
and
\begin{equation}\label{Xs}
\hat\P_s^\top\Omega\hat\I_s = \pmatrix{cccc}
\xi_0 &-\xi_1\\ \xi_1 &0&\ddots\\ &\ddots &\ddots &-\xi_{s-1}\\ &&\xi_{s-1} &0\endpmatrix =: X_s, \qquad \xi_i = \left(2\sqrt{|4i^2-1|}\right)^{-1},~i=0,\dots,s-1,
\end{equation} 
 that (\ref{gam}) can be rewritten as
\begin{equation}\label{gam1}
\bfgamma = \bfe_1\otimes p_0 + hX_s\otimes I\left(\bfvphi-\bfeta\right).
\end{equation}
Next, by plugging the right-hand side of (\ref{gam1}) in those of (\ref{eta})-(\ref{vfi}), and taking into account that
\begin{equation}\label{Ise1}
\hat\I_s\bfe_1 = \hat\bfc, \qquad \I_s\bfe_1 = \bfc,
\end{equation}
one obtains:
\begin{eqnarray}\label{eta1}
\bfeta &=& \P_s^\top \Omega\otimes I\,\nabla U\left( \bfuno_k\otimes q_0 + h\bfc\otimes p_0 + h^2\I_sX_s\otimes I (\bfvphi-\bfeta) \right),\\ 
\label{vfi1}
\bfvphi &=& \hat\P_s^\top \hat\Omega\otimes I\\ \nonumber
&& \left[\B\left( \bfuno_s\otimes q_0 + h\hat\bfc\otimes p_0 + h^2\hat\I_sX_s\otimes I (\bfvphi-\bfeta) \right)\cdot\left(\bfuno_s\otimes p_0 + h\hat\I_s\otimes I\left( \bfvphi -\bfeta\right)\right)\right].
\end{eqnarray}
Further, by defining the block vector 
\begin{equation}\label{psi}
\bfpsi \equiv \pmatrix{c} \psi_0\\ \vdots\\ \psi_{s-1}\endpmatrix := \bfvphi-\bfeta,
\end{equation}
subtracting (\ref{eta1}) from (\ref{vfi1}) provides us with the (block) vector equation
\begin{eqnarray}\nonumber
\bfpsi &=& \hat\P_s^\top \hat\Omega\otimes I \left[\B\left( \bfuno_s\otimes q_0 + h\hat\bfc\otimes p_0 + h^2\hat\I_sX_s\otimes I\, \bfpsi \right)\cdot\left(\bfuno_s\otimes p_0 + h\hat\I_s\otimes I\,\bfpsi\right)\right]\\ \label{psi1}
&&- \P_s^\top \Omega\otimes I\,\nabla U\left( \bfuno_k\otimes q_0 + h\bfc\otimes p_0 + h^2\I_sX_s\otimes I\, \bfpsi \right).
\end{eqnarray}
Once (\ref{psi1}) has been solved, one easily computes the first block entry of (\ref{gam1}) by taking into account (\ref{Xs}), so that the new approximations (\ref{newuh_1}) become
\begin{equation}\label{newuh1}
q_1 = q_0 +h p_0 +\frac{h^2}2\left(\psi_0 - \frac{1}{\sqrt{3}}\psi_1\right), \qquad p_1=p_0 +h\psi_0.
\end{equation}

\begin{rem} As is clear from (\ref{newuh1}), in order to obtain the new approximation $q_1$, one must have $s\ge2$, so that all LIM$(k,s)$ methods, with $k\ge s$, have order $2s\ge4$. Consequently, the methods here derived are completely different from the second-order method studied in \cite{LW2016}. 

\smallskip
Moreover, in the case where the magnetic field is zero (i.e.,  $\B=O$ in (\ref{psi1})), then the discrete problem (\ref{psi1})-(\ref{newuh1}) reduces to that generated by a HBVM$(k,s)$ method applied to the special second order problem ~$\ddot q +\nabla U(q)=0$,~ $q(0)=q_0$,~ $\dot q(0)=p_0$~ (see, e.g., \cite[Chapter\,4.1.2]{LIMbook2016}). 

\smallskip
Also, we observe that, introducing the notation
$$\pmatrix{c} x_1 \\ \vdots \\ x_s\endpmatrix \times \pmatrix{c} y_1 \\ \vdots \\ y_s\endpmatrix := \pmatrix{c} x_1\times y_1 \\ \vdots \\ x_s\times y_s\endpmatrix, \qquad  x_i,y_i\in\RR^3,\quad i=1,\dots,s,$$ then (\ref{psi1}) can be rewritten in a form closer to that of the original problem (\ref{cpd}), i.e.,
\begin{eqnarray}\nonumber
\bfpsi &=& \hat\P_s^\top \hat\Omega\otimes I \left[\left(\bfuno_s\otimes p_0 + h\hat\I_s\otimes I\,\bfpsi\right)\times 
L\left( \bfuno_s\otimes q_0 + h\hat\bfc\otimes p_0 + h^2\hat\I_sX_s\otimes I\, \bfpsi \right)\right]\\ \label{psi2}
&&- \P_s^\top \Omega\otimes I\,\nabla U\left( \bfuno_k\otimes q_0 + h\bfc\otimes p_0 + h^2\I_sX_s\otimes I\, \bfpsi \right).
\end{eqnarray}

\smallskip
Last, but not least, we stress that, in order to obtain an energy-conserving method of order $2s$, we need to solve, at each integration step,  the set of $3s$ algebraic equations (\ref{psi1}) (i.e., (\ref{psi2})), whose dimension is remarkably {\em independent of $k$}, as previously anticipated.
\end{rem}

Next result guarantees the existence and uniqueness of the solution of (\ref{psi1}).

\begin{theo}\label{hsmall}
Let $B(q)p$ and $\nabla U(q)$ in (\ref{cpd1}) be continuous and satisfy a Lipschitz condition with constant $\mu$ 
(with respect to $p$ and $q$) given by
$$
\left\| B(\bar q) \bar p - B( q) p \right\| \le \mu \left(\| \bar q - q \|  + \|\bar p -  p \| \right), \quad
\left\| \nabla U(\bar q) - \nabla U(q) \right\| \le \mu \| \bar q - q\|.
$$
Then, if the stepsize $h$ is small enough to satisfy
\begin{equation}\label{hmu}
h \mu\left[ \| \hat\I_s \| \,\| \hat\P_s^\top \hat\Omega_s \| + 
h\|X_s\| \left( \| \hat\I_s \| \, \| \hat\P_s^\top \hat\Omega_s \| +
\| \I_s \| \, \| \P_s^\top  \Omega_s \|
\right)\right] < 1,
\end{equation}
there exists a unique solution of (\ref{psi1}), and the fixed-point iteration 
\begin{eqnarray}\nonumber
\bfpsi^{\ell+1} &=& \hat\P_s^\top \hat\Omega\otimes I \left[\B\left( \bfuno_s\otimes q_0 + h\hat\bfc\otimes p_0 + h^2\hat\I_sX_s\otimes I\, \bfpsi^\ell \right)\cdot\left(\bfuno_s\otimes p_0 + h\hat\I_s\otimes I\,\bfpsi^\ell\right)\right]\\ \label{fpit}
&&- \P_s^\top \Omega\otimes I\,\nabla U\left( \bfuno_k\otimes q_0 + h\bfc\otimes p_0 + h^2\I_sX_s\otimes I\, \bfpsi^\ell \right), \qquad \ell=0,1,\dots,
\end{eqnarray}
converges to it. 
\end{theo}
\proof
First, note that the Lipschitz conditions ensure that the same functions with block vector arguments 
are also Lipchitz with the same constant. Then, it is straightforward to  prove that the function defined by the right-hand side of (\ref{psi1}) is a contraction with respect to $\bfpsi$ if $h$ satisfies (\ref{hmu}). Consequently, the Fixed-Point Theorem ensures the existence and uniqueness of the solution, as well as the convergence of the iteration (\ref{fpit}).\,\QED

\subsection{Symmetry}
The compact structure of the discrete problem (\ref{psi1})-(\ref{newuh1}) allows us to prove the important property of symmetry of the method, under the following symmetry condition of the abscissae,
\begin{equation}\label{cis}
\hat c_\ell = 1-\hat c_{s-\ell+1}, \quad \ell=1,\dots,s,\qquad\qquad c_\ell = 1-c_{k-\ell+1},\quad \ell=1,\dots,k, 
\end{equation}
which is clearly satisfied by the choice of the Gauss-Legendre abscissae $P_s(\hat c_\ell)=0$, $\ell=1,\dots,s$, $P_k(c_\ell) =0$, $\ell=1,\dots,k$. For interpolatory quadrature formulae, in turn, (\ref{cis}) implies the symmetry of the weights:
\begin{equation}\label{bis}
\hat b_\ell = \hat b_{s-\ell+1}, \quad \ell=1,\dots,s,\qquad\qquad b_\ell = b_{k-\ell+1},\quad \ell=1,\dots,k. 
\end{equation}
We also need to define the following matrices, 
\begin{equation}\label{PD}
\hat P = \pmatrix{ccc} & &1\\&\udots\\ 1\endpmatrix, \, \hat D = \pmatrix{ccc} (-1)^0\\ &\ddots\\ && (-1)^{s-1}\endpmatrix~\in~\RR^{s\times s},
 \quad P = \pmatrix{ccc} & &1\\&\udots\\ 1\endpmatrix~\in~\RR^{k\times k},
\end{equation}
and related preliminary results.

\begin{lem}\label{PIX}
With reference to the matrices defined in (\ref{PIO}), (\ref{Xs}), and (\ref{PD}), and the vector $\bfe_1$ defined in (\ref{e1}), one has:
\begin{eqnarray*} 
\hat D^2 ~=~\hat P^2 &=&I_s,\qquad P^2 ~=~I_k,\qquad \hat D X_s \hat D = X_s^\top \equiv \bfe_1\bfe_1^\top -X_s.\\[1mm]
\hat P\hat\I_s \hat D &=& \bfuno_s\bfe_1^\top -\hat\I_s, \qquad P\I_s \hat D ~=~ \bfuno_k\bfe_1^\top -\I_s,\\[1mm] 
\hat P\hat\P_s \hat D &=& \hat\P_s, \qquad\qquad\quad P\P_s \hat D ~=~ \P_s,\\[1mm] 
\hat P\hat\Omega\hat P &=&\hat\Omega, \qquad\qquad\quad\, P\Omega P ~~~=~\Omega.
\end{eqnarray*}
\end{lem}
\proof The properties on the first line follows by the fact that $\hat D,\hat P,P$ are symmetric and orthogonal, and from (\ref{Xs}). The properties on the subsequent two lines derive from the following symmetries of the Legendre polynomials: $$P_j(1-c) = (-1)^jP_j(c), \qquad \int_0^{1-c}P_j(x)\dd x = \delta_{j0}-(-1)^j\int_0^c P_j(x)\dd x, \qquad j=0,1,\dots.$$ At last, the properties on the last line follow from (\ref{bis}).\,\QED\bigskip

We recall that the method (\ref{psi1})-(\ref{newuh1}) is symmetric if, when starting from $(q_1,p_1)$ with timestep $-h$, it brings back to $(q_0,p_0)$.\footnote{The proof of symmetry will be akin  to that done for HBVMs in  \cite[Theorem\,3.11]{LIMbook2016}.}

\begin{theo}\label{symth}
Under the hypotheses (\ref{cis}), the method (\ref{psi1})-(\ref{newuh1}) is symmetric.
\end{theo}
\proof
By using the method (\ref{psi1})-(\ref{newuh1}) for solving (\ref{cpd1})-(\ref{Bq}) starting from $(q_1,p_1)$ with stepsize $-h$, we obtain the equations

\begin{eqnarray}\nonumber
\bar\bfpsi \equiv \pmatrix{c}\bar\psi_0\\ \vdots\\ \bar\psi_{s-1}\endpmatrix&=& \hat\P_s^\top \hat\Omega\otimes I \left[\B\left( \bfuno_s\otimes q_1 - h\hat\bfc\otimes p_1 + h^2\hat\I_sX_s\otimes I\, \bar\bfpsi \right)\cdot\left(\bfuno_s\otimes p_1 - h\hat\I_s\otimes I\,\bar\bfpsi\right)\right]\\  \label{psi3}
&&- \P_s^\top \Omega\otimes I\,\nabla U\left( \bfuno_k\otimes q_1 - h\bfc\otimes p_1 + h^2\I_sX_s\otimes I\, \bar\bfpsi \right),\\ \nonumber
\bar q_0 &=& q_1 -h p_1 +\frac{h^2}2\left(\bar\psi_0 - \frac{1}{\sqrt{3}}\bar\psi_1\right), \qquad \bar p_0~=~p_1 -h\bar\psi_0.
\end{eqnarray}
We have then to prove that
\begin{equation}\label{finesym}
\bar q_0 = q_0, \qquad\quad \bar p_0 = p_0.
\end{equation}
Preliminarily, we observe that the symmetry conditions (\ref{cis}) can be respectively rewritten  as $$\hat P\hat\bfc =\bfuno_s-\hat\bfc, \qquad P\bfc = \bfuno_k-\bfc.$$ Moreover, from the last two equations in (\ref{psi3}) and (\ref{newuh1}), one has
\begin{equation}\label{newh3}
\bar q_0 = q_0 -\frac{h^2}2\left[ \left(\psi_0-\bar\psi_0\right) + \frac{1}{\sqrt{3}}\left( \psi_1+\bar\psi_1\right)\right],\qquad \bar p_0 = p_0 + h\left( \psi_0-\bar\psi_0\right).
\end{equation}
Then, (\ref{finesym}) follows if we show that (see (\ref{psi}))
\begin{equation}\label{psistar}
\bar\psi_j=(-1)^j\psi_j =:\psi_j^*, \qquad j=0,\dots,s-1,\qquad \Leftrightarrow\qquad \bar\bfpsi = \hat D\otimes I\bfpsi =:\bfpsi^* \equiv \pmatrix{c} \psi_0^*\\ \vdots\\ \psi_{s-1}^*\endpmatrix.
\end{equation}
We shall prove this statement by verifying that  the (block) vector $\bfpsi^*$ satisfies the very same equation (\ref{psi3}) which defines $\bar\bfpsi$. From (\ref{psi1}) and Lemma~\ref{PIX}, one has:
\begin{eqnarray*}
\bfpsi^*&=& \hat D\otimes I\, \bfpsi\\
&=&\hat D\, \hat\P_s^\top \hat\Omega\otimes I \left[\B\left( \bfuno_s\otimes q_0 + h\hat\bfc\otimes p_0 + h^2\hat\I_sX_s\hat D\otimes I\, \bfpsi^* \right)\cdot\left(\bfuno_s\otimes p_0 + h\hat\I_s\hat D\otimes I\,\bfpsi^*\right)\right]\\
&&- \hat D\,\P_s^\top \Omega\otimes I\,\nabla U\left( \bfuno_k\otimes q_0 +h\bfc\otimes p_0 + h^2\I_sX_s \hat D \otimes I\, \bfpsi^* \right)\\
&=& \hat\P_s^\top \hat\Omega \hat P\otimes I \left[\B\left( \bfuno_s\otimes q_0 + h\hat\bfc\otimes p_0 + h^2\hat\I_sX_s\hat D\otimes I\, \bfpsi^* \right)\cdot\left(\bfuno_s\otimes p_0 + h\hat\I_s\hat D\otimes I\,\bfpsi^*\right)\right]\\
&&- \P_s^\top \Omega P\otimes I\,\nabla U\left( \bfuno_k\otimes q_0 + h\bfc\otimes p_0 + h^2\I_sX_s\hat D\otimes I\, \bar\bfpsi \right)\\
&=& \hat\P_s^\top \hat\Omega \otimes I \left[\B\left( \hat P\bfuno_s\otimes q_0 + h\hat P\hat\bfc\otimes p_0 + h^2\hat P\hat\I_s\hat D^2X_s\hat D\otimes I\, \bfpsi^* \right)\cdot \left(\hat P\bfuno_s\otimes p_0 + h\hat P\hat\I_s\hat D\otimes I\,\bfpsi^*\right)\right]\\
&&- \P_s^\top \Omega \otimes I\,\nabla U\left( P\bfuno_k\otimes q_0 + hP\bfc\otimes p_0 + h^2P\I_s\hat D^2X_s\hat D\otimes I\, \bfpsi^* \right)\\
&=& \hat\P_s^\top \hat\Omega \otimes I \left[\B\left( \bfuno_s\otimes q_0 + h(\bfuno_s-\hat\bfc)\otimes p_0 + h^2\left(\bfuno_s\bfe_1^\top-\hat\I_s\right) X_s^\top\otimes I\, \bfpsi^* \right) \cdot \right.\\
&&\left.  \left(\bfuno_s\otimes p_0 + h\left(\bfuno_s\bfe_1^\top-\hat\I_s\right) \otimes I\,\bfpsi^*\right)\right]\\
&&- \P_s^\top \Omega \otimes I\,\nabla U\left( \bfuno_k\otimes q_0 + h(\bfuno_k-\bfc)\otimes p_0 +h^2\left(\bfuno_k\bfe_1^\top- \I_s\right) X_s^\top\otimes I\, \bfpsi^* \right) ~=~(*).\\
\end{eqnarray*}
Next, upon observing that (see (\ref{Xs}), (\ref{newuh1}), and (\ref{psistar}))
\begin{eqnarray*}
&&q_0+hp_0+h^2\bfe_1^\top X_s^\top\otimes I\,\bfpsi^* = q_0+hp_0+\frac{h^2}2\left(\psi_0^*+\frac{1}{\sqrt{3}}\psi_1^*\right) = q_0+hp_0+\frac{h^2}2\left(\psi_0-\frac{1}{\sqrt{3}}\psi_1\right) = q_1,\\
&& p_0+h\bfe_1^\top\otimes I\, \bfpsi^* = p_0+h\psi_0^* =p_0+h\psi_0 = p_1,
\end{eqnarray*}
we obtain, taking into account again that $p_0=p_1-h\psi_0^*$,
\begin{eqnarray*}
(*) &=& \hat\P_s^\top \hat\Omega \otimes I \left[\B\left( \bfuno_s\otimes q_1 - h\hat\bfc\otimes p_1  + h^2\left[\hat\bfc\bfe_1^\top -\hat\I_s X_s^\top\right]\otimes I\, \bfpsi^* \right)\cdot\left(\bfuno_s\otimes p_1 - h\hat\I_s\otimes I\,\bfpsi^*\right)\right]\\
&&- \P_s^\top \Omega \otimes I\,\nabla U\left( \bfuno_k\otimes q_1 - h\bfc\otimes p_1 + h^2\left[\bfc\bfe_1^\top -\I_s X_s^\top\right]\otimes I\, \bfpsi^* \right).
\end{eqnarray*}
The statement then follows by considering that (see (\ref{PIO}), (\ref{Xs}),  (\ref{Ise1}), and Lemma~\ref{PIX})
$$\left[\hat\bfc\bfe_1^\top -\hat\I_sX_s^\top\right] = \hat\I_s\left( \bfe_1\bfe_1^\top -\bfe_1\bfe_1^\top+X_s\right) = \hat\I_s X_s$$
and, similarly, ~$\left[\bfc\bfe_1^\top -\I_s X_s^\top\right] = \I_s X_s.\,\QED$\medskip

\subsection{Solving the discrete problem}\label{blesec}

When the stepsize $h$ satisfies the conditions of Theorem~\ref{hsmall}, we proved that the fixed-point iteration (\ref{fpit}) converges to the solution of (\ref{psi1}). However, sometimes a simplified Newton iteration may be more appropriate for solving the equation
\begin{eqnarray}\nonumber 
\lefteqn{F(\bfpsi) ~:=}\\ \nonumber 
&&\bfpsi - \hat\P_s^\top \hat\Omega\otimes I \left[\B\left( \bfuno_s\otimes q_0 + h\hat\bfc\otimes p_0 + h^2\hat\I_sX_s\otimes I\, \bfpsi \right)\cdot \left(\bfuno_s\otimes p_0 + h\hat\I_s\otimes I\,\bfpsi\right)\right]\\ \label{Fpsi}
&&+\, \P_s^\top \Omega\otimes I\,\nabla U\left( \bfuno_k\otimes q_0 + h\bfc\otimes p_0 + h^2\I_sX_s\otimes I\, \bfpsi \right) ~=~\bfzero.
\end{eqnarray}
Neglecting the $O(h^2)$ terms in the simplified Jacobian, and taking into account (\ref{Xs}), one then obtains the iteration
\begin{equation}\label{sNit}
\mbox{solve:~}\left[I_s\otimes I-hX_s\otimes B(q_0)\right] \Delta\bfpsi^\ell = -F(\bfpsi^\ell), \qquad \ell=0,1,\dots.
\end{equation}
This iteration, though straightforward, requires the factorization of a matrix having dimension $3s$. This complexity can be reduced, e.g., by using variants akin to those used for HBVMs \cite{BIT2011,BFCI2014}.   In particular, we consider the {\em blended iteration} associated to (\ref{sNit}), which, by setting\,\footnote{We refer to \cite{BM2002} for a more comprehensive analysis of {\em blended methods}.} 
\begin{equation}\label{teta}
\Theta = (I-h\rho_s B(q_0))^{-1}, \qquad \rho_s = \min_{\lambda\in\sigma(X_s)}|\lambda|,
\end{equation}
 reads    
\begin{equation}\label{blit}
\bfb^\ell := -F(\bfpsi^\ell),\quad \bfb_1^\ell := \rho_s X_s^{-1}\otimes I\,\bfb^\ell,\quad
\Delta\bfpsi^\ell = I_s\otimes\Theta\left[\bfb_1^\ell + I_s\otimes\Theta\left(\bfb^\ell-\bfb_1^\ell\right)\right], \quad \ell=0,1,\dots.
\end{equation}
In so doing, only the $3\times 3$ matrix $\Theta$ in (\ref{teta}) needs to be computed.

Sometimes, the electric field may be much stronger than the magnetic one. In such a case, the simplified Newton iteration (\ref{sNit}) becomes, by neglecting the contribution of the magnetic field, and considering that $\P_s^\top\Omega\I_s=X_s$, as defined in (\ref{Xs}),
\begin{equation}\label{sNit1}
\mbox{solve:~}\left[I_s\otimes I+h^2X_s^2\otimes \nabla^2 U(q_0)\right] \Delta\bfpsi^\ell = -F(\bfpsi^\ell), \qquad \ell=0,1,\dots.
\end{equation}
We can again use a corresponding {\em blended iteration} \cite{LIMbook2016,BM2007},
\begin{equation}\label{blit1}
\bfb^\ell := -F(\bfpsi^\ell),\quad \bfb_1^\ell := \rho_s^2 X_s^{-2}\otimes I\,\bfb^\ell,\quad
\Delta\bfpsi^\ell = I_s\otimes\Theta_1\left[\bfb_1^\ell + I_s\otimes\Theta_1\left(\bfb^\ell-\bfb_1^\ell\right)\right], \quad \ell=0,1,\dots,
\end{equation}
where $\rho_s$ is the same parameter defined in (\ref{teta}), and 
\begin{equation}\label{teta1}
\Theta_1 = (I+h^2\rho_s^2\nabla^2U(q_0))^{-1}. 
\end{equation}
Consequently, also in this case, only a $3\times 3$ matrix need to be factored.

\begin{rem} We observe that, in the case of the more general problem (\ref{general}), the previous arguments remains formally the same, with the only difference that, now, the matrices $\Theta$ and $\Theta_1$ defined in (\ref{teta}) and (\ref{teta1}), respectively, have dimension $m\times m$.\end{rem}

\begin{rem}[The case of a constant $B$]
When the magnetic field $L(q)$ is uniform, then with reference to matrix $B(q)$ defined in (\ref{Bq}), one has:
$$B(q)\equiv B=-B^\top.$$ In this case, taking into account (\ref{e1})-(\ref{Xs}), one has that the discrete problem (\ref{psi1}) simplifies to
\begin{equation}\label{psi1_1}
\bfpsi = \bfe_1\otimes Bp_0 + hX_s\otimes B\,\/\bfpsi - \P_s^\top \Omega\otimes I\,\nabla U\left( \bfuno_k\otimes q_0 + h\bfc\otimes p_0 + h^2\I_sX_s\otimes I\, \bfpsi \right).
\end{equation}
Moreover, also the blended iteration (\ref{teta})-(\ref{blit}) greatly simplifies, since the matrix $\Theta$ turns out to be given by  ~$\Theta=(I-h\rho_sB)^{-1}$~ and, therefore, it is constant for all time-steps. 
\end{rem}

\section{Numerical tests}\label{num}
We here report a few numerical tests, aimed at assessing the theoretical findings of the previous sections. We also compare the LIM$(k,s)$ method (\ref{psi1})-(\ref{newuh1}) with the Boris method.\footnote{An efficient implementation of this latter method is explained in \cite{HL2018}.} 
All numerical tests have been done on a 2.8GHz Intel core i7 computer with 16GB of memory, running Matlab 2017b.

\begin{figure}[t]
\centerline{\includegraphics[width=7cm,height=6cm]{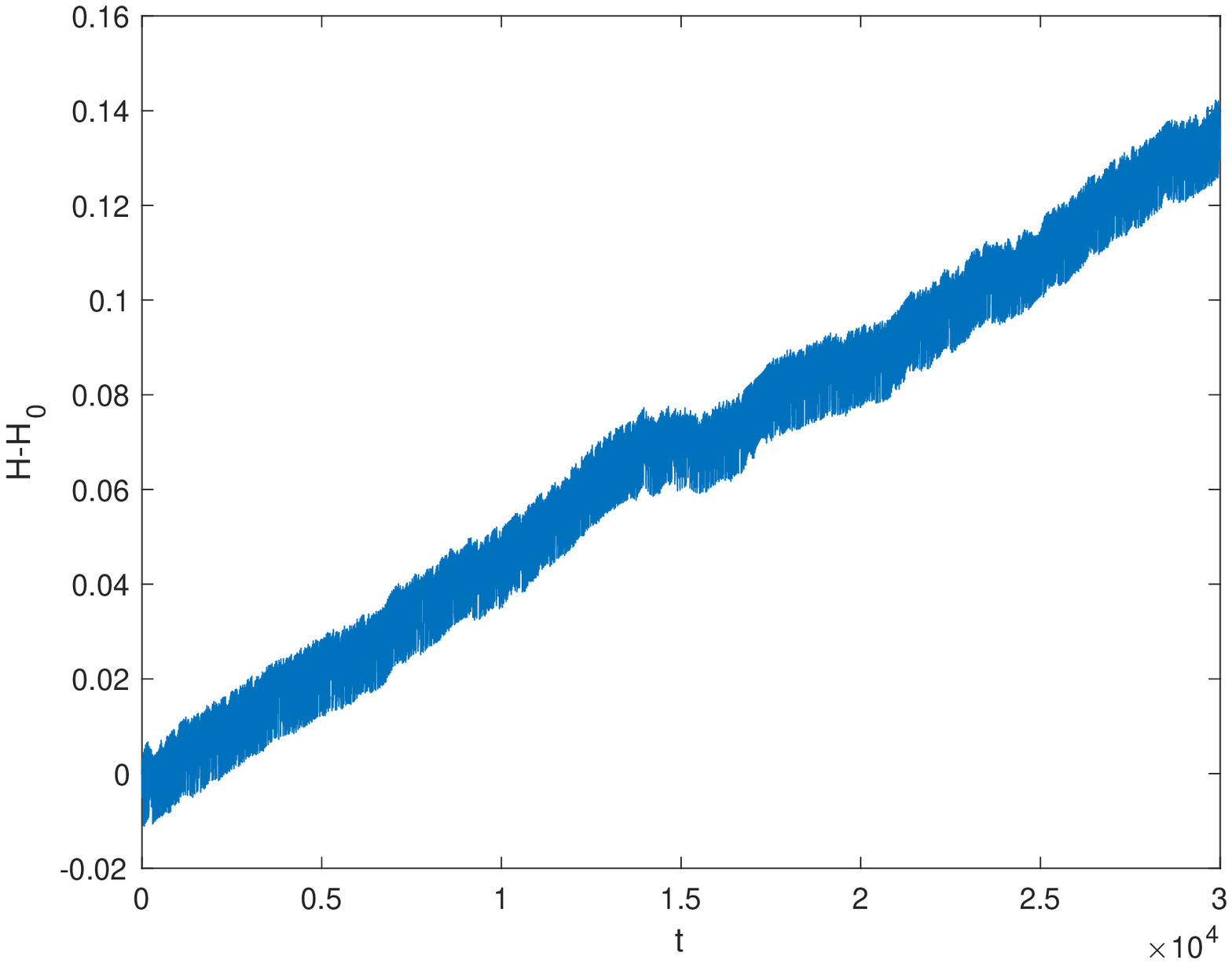}\qquad\includegraphics[width=7cm,height=6.2cm]{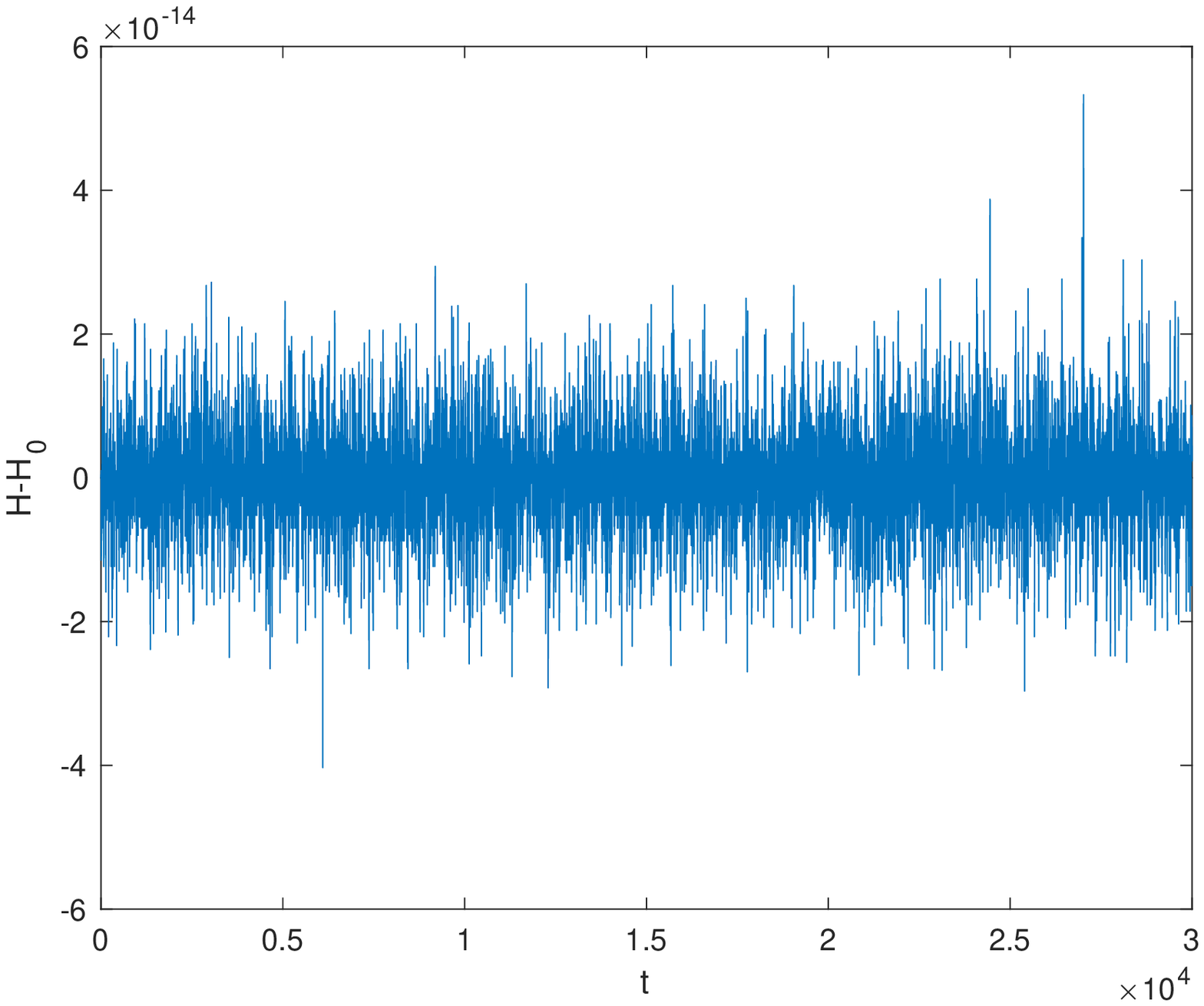}}
\caption{Hamiltonian error when solving problem (\ref{U1})--(\ref{L1}) with stepsize $h=10^{-2}$. Left plot: Boris method. Right plot: LIM(4,2).}\label{fig1}

\bigskip
\bigskip

\centerline{\includegraphics[width=12cm,height=9cm]{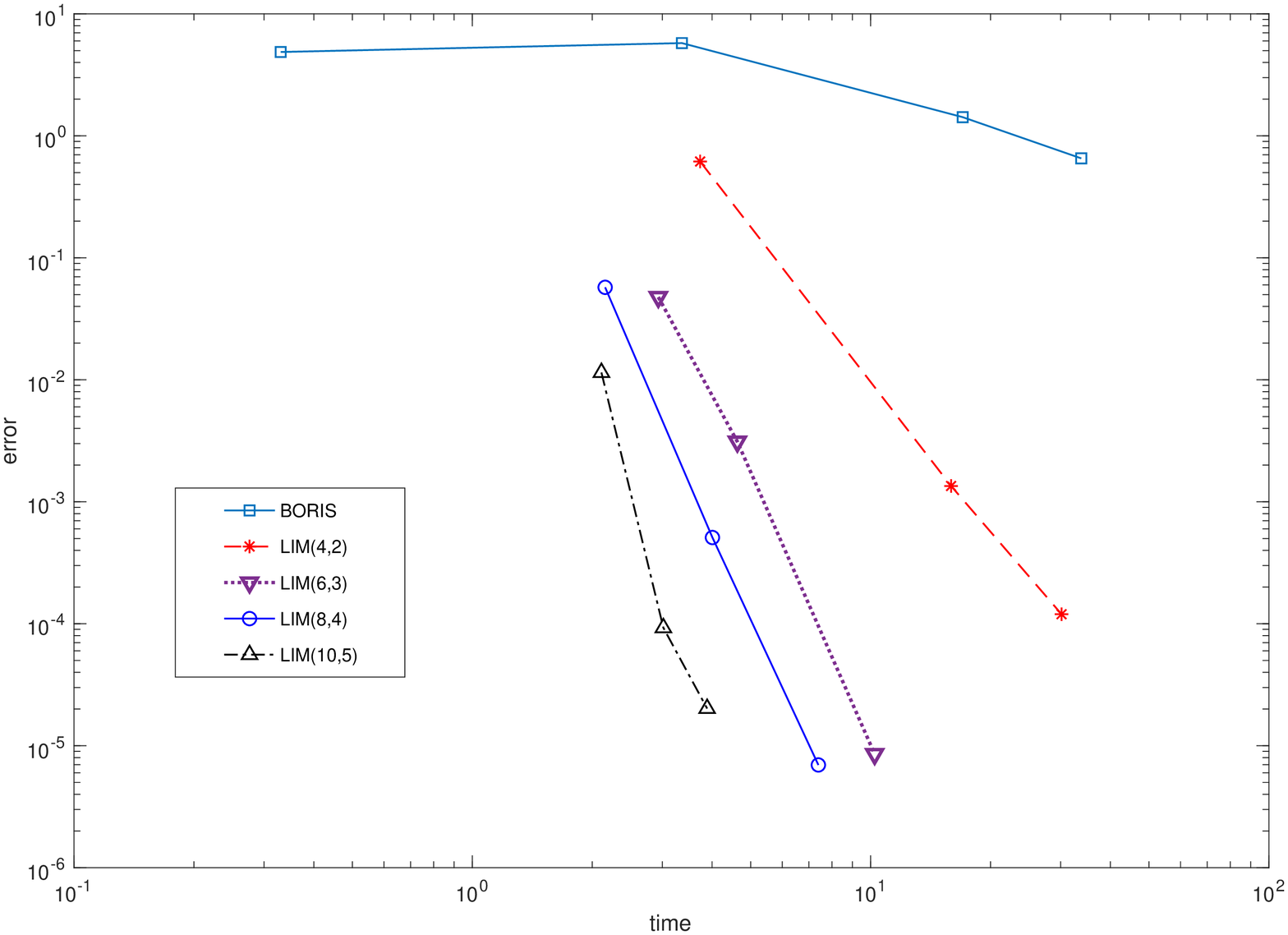}}
\caption{Solution error vs. execution time for problem (\ref{U1})--(\ref{L1}) at $t=100$.}\label{fig1_1}
\end{figure}

\paragraph*{Example 1.}
We start considering the problem (\ref{cpd})--(\ref{Bq}) with:\footnote{Hereafter, $q_i$ will denote the $i$th entry of $q$.}
\begin{eqnarray}\label{U1}
U(q) &=& q_1^3-q_2^3+\frac{1}5q_1^4+q_2^4+q_3^4,\\  \label{q0p0}
q(0) &=& \pmatrix{ccc} 0, & 1, & 0.1\endpmatrix^\top,\qquad
p(0) ~=~ \pmatrix{ccc} 0.09, & 0.55, & 0.3\endpmatrix^\top,\\ \label{L1}
L(q) &=& \pmatrix{ccc} 0, & 0,& \sqrt{q_1^2+q_2^2}\endpmatrix^\top,
\end{eqnarray} 
for which  in \cite{HL2018} it has been shown that the Boris algorithm exhibits a $O(th^2)$ drift in the energy. For this problem, the LIM$(2s,s)$ method turns out to be energy-conserving, according to the result of Theorem~\ref{newH1H0}, since $U(q)$ is a polynomial of degree 4. The corresponding order is  $2s$, as stated in Theorem~\ref{newy1yh}. The drift for the Boris method is confirmed by the left-plot in Figure~\ref{fig1}, where the Hamiltonian error for such method, using a stepsize $h=10^{-2}$ over the interval $[0,3\cdot 10^4]$, is shown:  the obtained plot, upon scaling by $h^2$ and reversing the sign, perfectly fits those in \cite[Fig.\,1]{HL2018}, thus confirming the $O(th^2)$ energy drift.
In the right-plot of the same figure, there is the Hamiltonian error for the LIM(4,2) method using the same stepsize: as one may see, in such a case  energy-conservation is gained (clearly, up to round-off errors). It must be observed, however, that the execution time for the Boris method is about 11 times smaller than that for the LIM(4,2) method (101 sec vs. 1149 sec). Nevertheless, the availability of arbitrarily high-order energy-conserving line-integral methods (i.e., LIM$(2s,s)$, $s=2,3,\dots$), coupled with the efficient nonlinear iteration (\ref{teta})-(\ref{blit}) studied in Section~\ref{blesec}, allows to recover on the efficiency of the methods, when we compare the solution error vs. the execution time. For this, we fix a shorter interval, say $[0,100]$, and we compare various energy-conserving line integral methods, along with the Boris method, in terms of solution error vs. execution time. The obtained results are plotted in Figure~\ref{fig1_1}, from which one realizes that the higher order LIMs are more effective than the lower order ones and, in any case, all of them perform better than the Boris method, which is the lowest order method, among those considered.

\paragraph*{Example 2.}
Next, we consider the problem defined by (\ref{U1})-(\ref{q0p0}), and
\begin{equation}\label{L2}
L(q) = \frac{1}2\pmatrix{ccc} q_2-q_3,& q_1+q_3,& q_2-q_1\endpmatrix^\top,
\end{equation}
for which the energy error of the Boris algorithm defines a random walk \cite{HL2018}. For this problem, since $U(q)\in\Pi_4$ and $L(q)\in\Pi_1$, one has that the quadratures (\ref{hg})--(\ref{he}) are exact for all LIM$(2s,s)$ methods, so that the discrete method (\ref{newuv})-(\ref{newuvh}) coincides with the approximation procedure (\ref{uv})-(\ref{uvh}). In particular, all of such methods are energy-conserving and have order $2s$. In Table~\ref{tab2} we list the obtained results when solving the problem on the interval $[0,25]$ with stepsize $h=0.05/n$: as one may see, the expected convergence order is confirmed.

\begin{table}[t]
\caption{Maximum absolute errors when solving problem (\ref{U1})-(\ref{q0p0}) and (\ref{L2}) on the interval $[0,25]$ with stepsize $h=0.05/n$:
$e_y$ = solution error; $e_H$ = Hamiltonian error.}\label{tab2}

\smallskip
\centerline{\begin{tabular}{|r|rcrc|rcr|rcr|}
\hline
        &\multicolumn{4}{c|}{Boris} &\multicolumn{3}{c|}{LIM(4,2)} &\multicolumn{3}{c|}{LIM(6,3)}\\
        \hline
$n$ & $e_y$ & rate & $e_H$ & rate & $e_y$ & rate & $e_H$ & $e_y$ & rate & $e_H$ \\ 
\hline
 1 & 3.30e\,00 & --- & 1.82e-01 & --- & 1.86e-02 & --- & 2.25e-14 & 1.81e-05 & --- & 2.14e-14 \\ 
 2 & 8.67e-01 & 1.9 & 4.53e-02 & 2.0 & 1.17e-03 & 4.0 & 3.03e-14 & 2.84e-07 & 6.0 & 2.30e-14 \\ 
 4 & 2.18e-01 & 2.0 & 1.13e-02 & 2.0 & 7.30e-05 & 4.0 & 2.03e-14 & 4.10e-09 & 6.1 & 3.12e-14 \\ 
 8 & 5.46e-02 & 2.0 & 2.82e-03 & 2.0 & 4.56e-06 & 4.0 & 1.81e-14 & 5.53e-10 & *** & 2.68e-14 \\ 
16 & 1.37e-02 & 2.0 & 7.05e-04 & 2.0 & 2.85e-07 & 4.0 & 1.94e-14 & 5.27e-10 & *** & 2.83e-14 \\ 
\hline
\end{tabular}}
\end{table}

\begin{figure}[t]
\centerline{\includegraphics[width=7cm,height=6.5cm]{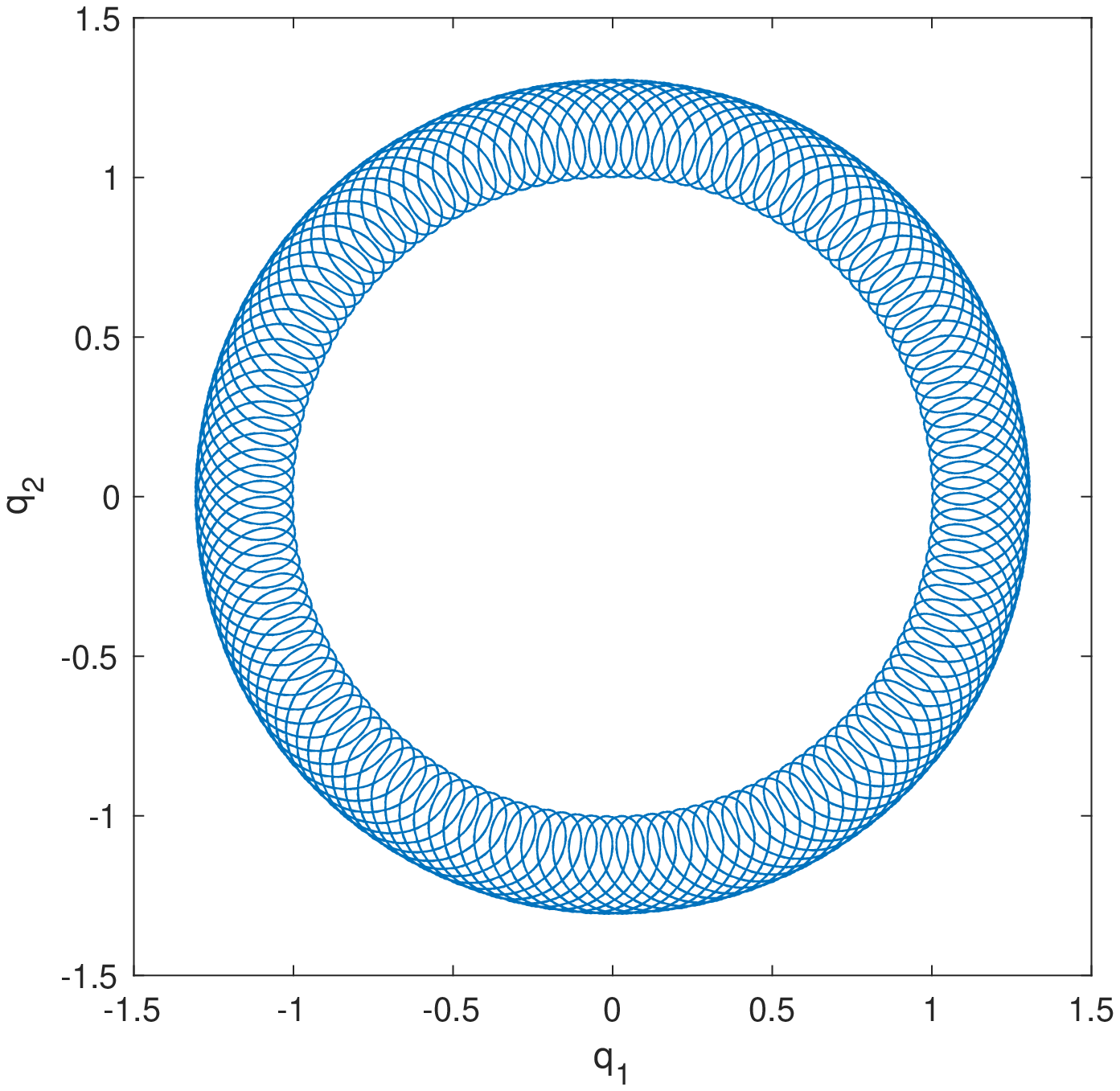}\qquad\includegraphics[width=7cm,height=6.5cm]{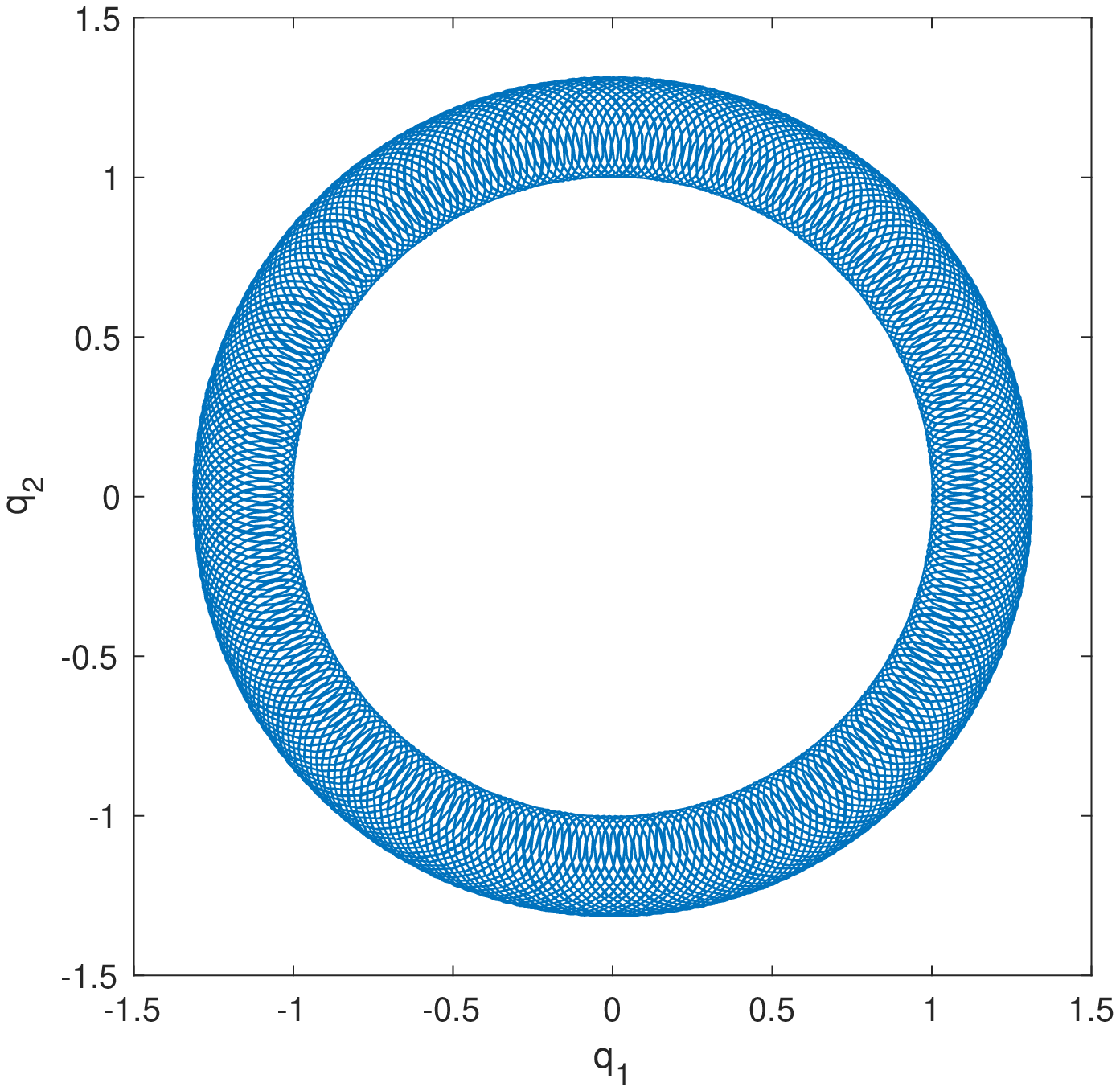}}

\bigskip
\bigskip

\centerline{\includegraphics[width=7cm,height=6.5cm]{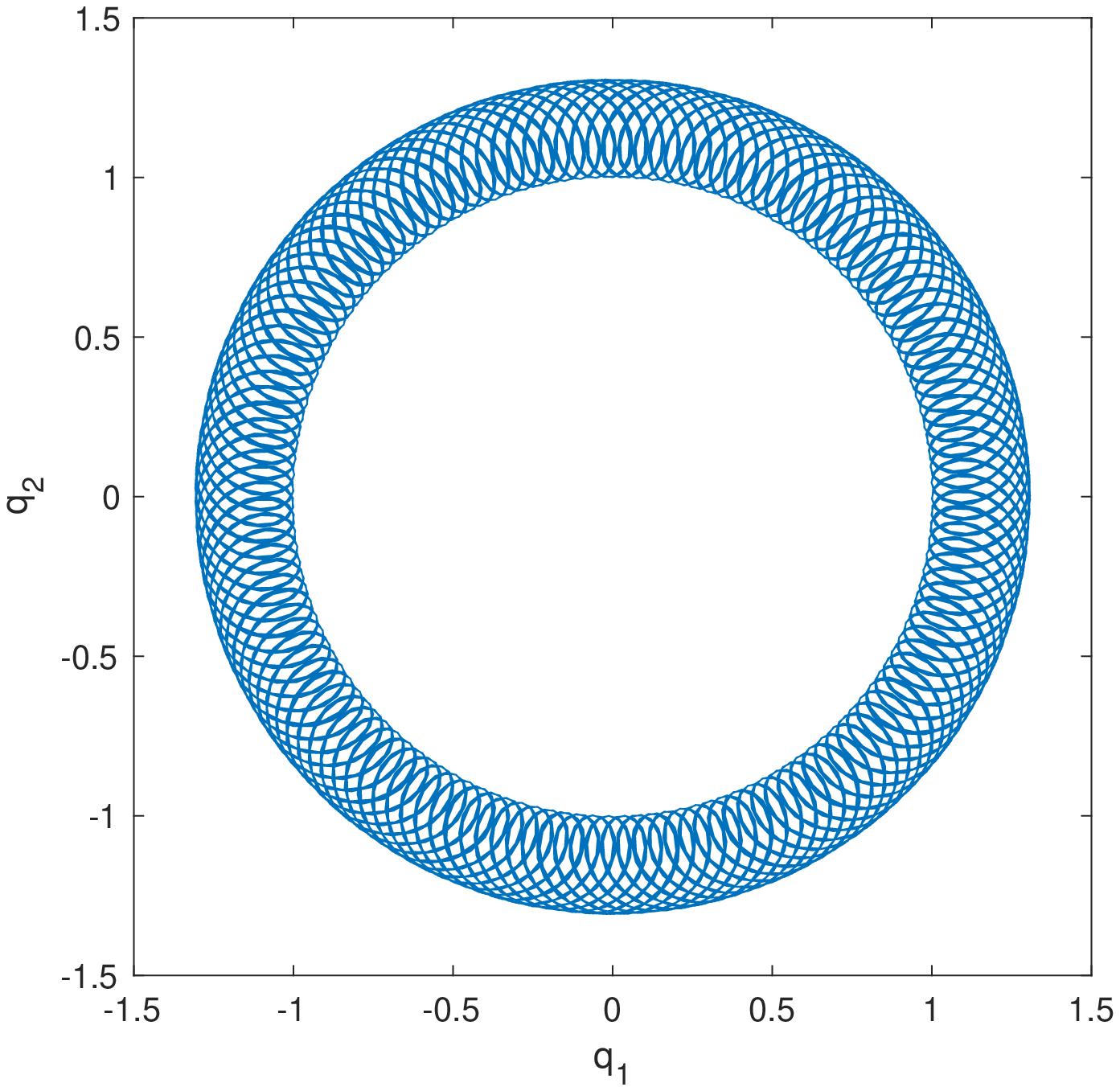}\qquad\includegraphics[width=7cm,height=6.5cm]{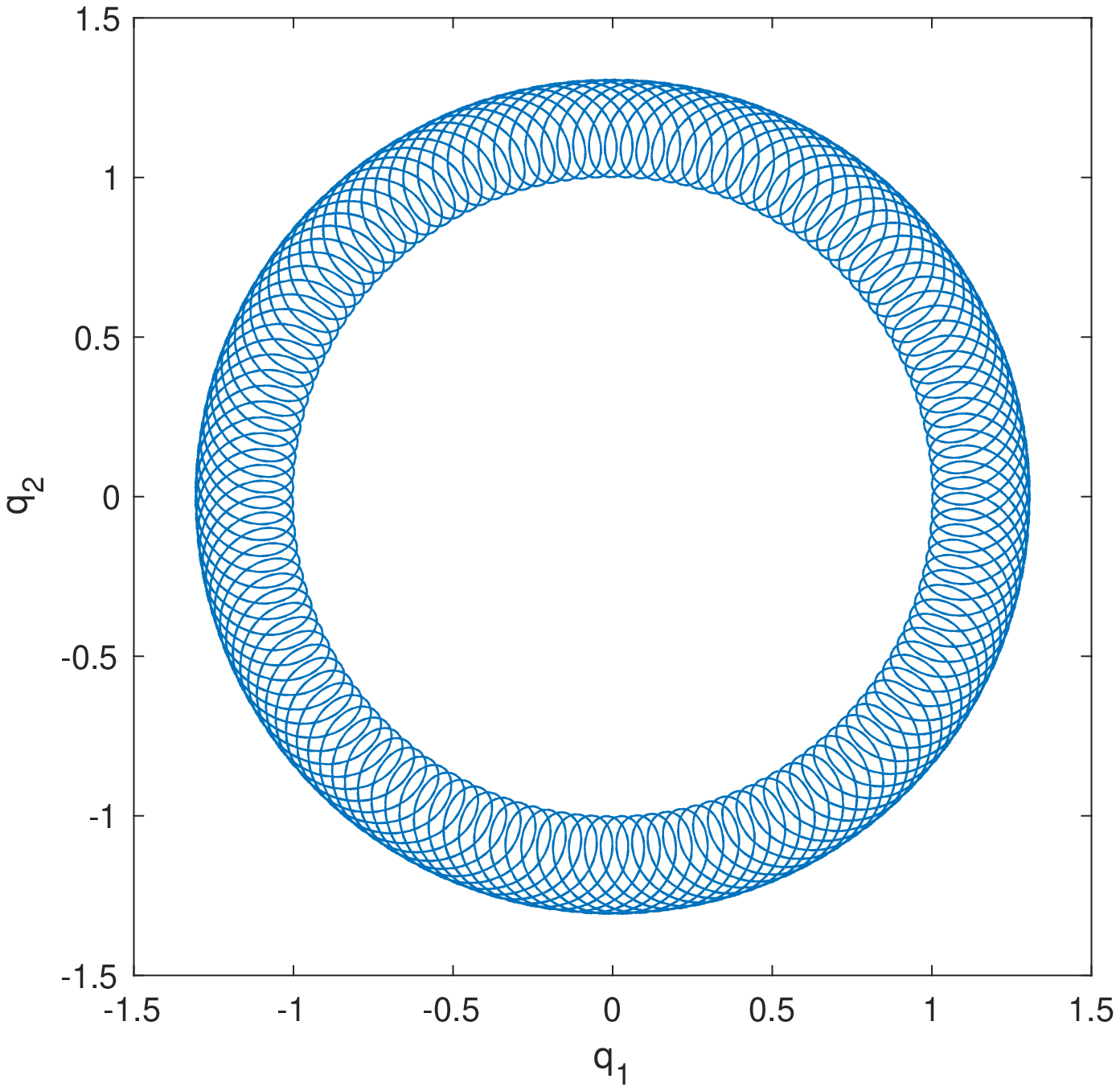}}
\caption{Problem (\ref{U-2D})-(\ref{q0p0-2D}) and (\ref{L1}) solved on the interval $[0,10^3\pi]$ by using a stepsize $h=\pi/10$.
Upper plots: reference solution (left) and Boris method (right). Lower plots: LIM(4,2) (left) and LIM(6,3) (right).}\label{fig3}
\end{figure}

\paragraph*{Example 3.}
Finally, we consider the 2D dynamics of a charged particle in a static, non-uniform electromagnetic field. The model is an important application in the study of the single particle motion and the guiding center dynamics \cite{LW2016,QZXLST2013}. The problem that we consider is defined by:\begin{eqnarray}\label{U-2D}
U(q) &=& \left[10\left(q_1^2+q_2^2\right)\right]^{-1},\\  \label{q0p0-2D}
q(0) &=& \pmatrix{ccc} 0, & 1, & 0\endpmatrix^\top,\qquad
p(0) ~=~ \pmatrix{ccc} 0.1, & 0.01, & 0\endpmatrix^\top, 
\end{eqnarray} 
and $L(q)$ defined as in (\ref{L1}). In this case, the motion, depicted in the upper-left plot of Figure~\ref{fig3} for $t\in[0,10^3\pi]$, occurs in the $(q_1,q_2)$-plane,  and another invariant is given by the angular momentum \cite{HSLQ2015},
$$
M(q,p) = q_1p_2-q_2p_1-\frac{(q_1^2+q_2^2)^{\frac{3}2}}3.
$$
We solve numerically this problem  by using the Boris and the LIM$(2s,s)$ methods, $s=2,3,4,5$, with stepsize $h=\pi/10$. 
The numerical solution of the Boris method is depicted in the upper-right plot of Figure~\ref{fig3}, whereas in the lower plots one finds the numerical solutions computed by the LIM(4,2) method (left) and LIM(6,3) method (right): as is clear, the higher the order of the method, the better the numerical solution. For completeness, in Table~\ref{tab3} we also list the maximum solution, Hamiltonian, and angular momentum errors for all methods, along with the measured execution times, again confirming that the higher order methods are more effective. In particular, it is worth mentioning that all considered LIMs are practically energy-conserving and, remarkably, their execution times are almost the same.

\begin{table}[t]
\caption{Maximum absolute errors when solving problem (\ref{U-2D})-(\ref{q0p0-2D}) and (\ref{L1}) on the interval $[0,10^3\pi]$ with stepsize $h=\pi/10$: $e_y$ = solution error; $e_H$ = Hamiltonian error; $e_M$ = angular momentum error. The execution times are in sec.}\label{tab3}

\smallskip
\centerline{\begin{tabular}{|l|r|r|r|c|}
\hline
method & $e_y\qquad$ & $e_H\qquad$  & $e_M\qquad$& execution time \\ 
\hline
Boris        & 2.5611e\,00 & 1.1461e-03 & 1.5532e-02 &0.4\\
LIM(4,2)   & 2.4553e-02 & 4.1633e-17  & 3.5917e-07 &6.0\\
LIM(6,3)   & 3.2533e-05 & 4.1633e-17  & 8.4765e-10 &6.2\\
LIM(8,4)   & 3.4584e-08 & 4.1633e-17  & 1.8433e-12 &6.3\\
LIM(10,5) & 7.9031e-09 & 4.1633e-17  & 1.9790e-11 &6.4\\
\hline
\end{tabular}}
\end{table}

\section{Conclusions}\label{fine}
In this paper we have derived symmetric and arbitrarily high-order energy-conserving methods for the numerical simulation of the dynamics of a charged particle, within the framework of {\em Line Integral Methods}. A complete analysis of the methods, including their efficient implementation, has been given.  Possible extensions to more general problems have also been sketched. Numerical results  on significant test problems, including the guiding center dynamics, duly confirm the theoretical findings.

\subsection*{Acknowledgements} This work has been finalized during the visit of the first author at the I.U.M.A., Universidad de Zaragoza, whose financial support is acknowledged (MINECO project MTM2016-77735-C3-1-P). The first author also acknowledges the interesting discussions with Ernst Hairer, during the RSME meeting in Santander, concerning the Boris method.


\begin{thebibliography}{10}

\bibitem{ABI2015} P.\,Amodio, L.\,Brugnano, F.\,Iavernaro. Energy-conserving methods for Hamiltonian Boundary Value Problems and applications in astrodynamics. {\em Adv. Comput. Math.} {\bf 41} (2015) 881--905. 

\bibitem{Be2008} P.M.\,Bellan. {\em Fundamentals of Plasma Physics}. Cambridge University Press, 2008.

\bibitem{B1970} J.P.\,Boris. Relativistic plasma simulation-optimization of a hybrid code. {\em Proceeding of Fourth Conference on Numerical Simulations of Plasmas}, pages 3--67, November 1970.

\bibitem{BCMR2012} L.\,Brugnano, M.\,Calvo, J.I.\,Montijano, L.\,R\'andez. Energy preserving methods for Poisson systems. {\em J. Comput. Appl. Math.} {\bf 236} (2012) 3890--3904.

\bibitem{BFCI2014}  L.\,Brugnano, G.\,Frasca Caccia, F.\,Iavernaro. Efficient  implementation of Gauss collocation and Hamiltonian Boundary Value Methods. {\em Numer. Algorithms} {\bf 65} (2014) 633--650.

\bibitem{BGI2018} L.\,Brugnano, G.\,Gurioli, F.\,Iavernaro. Analysis of Energy and QUadratic Invariant Preserving (EQUIP) methods.  {\em J. Comput. Appl. Math.} {\bf 335} (2018) 51--73.

\bibitem{BGIW2018} L.\,Brugnano, G.\,Gurioli, F.\,Iavernaro, E.B.\,Weinm\"uller. Line integral solution of Hamiltonian systems with holonomic constraints. {\em Appl. Numer. Math.} {\bf 127} (2018) 56--77.

\bibitem{BI2012} L.\,Brugnano, F.\,Iavernaro.  Line Integral Methods which preserve all invariants of conservative problems.  {\em J. Comput. Appl. Math.} {\bf 236} (2012) 3905--3919.

\bibitem{LIMbook2016} L.\,Brugnano, F.\,Iavernaro. {\em Line Integral Methods for Conservative Problems}.  Chapman and Hall/CRC, Boca Raton, FL, 2016.

\bibitem{BI2018} L.\,Brugnano, F.\,Iavernaro. Line Integral Solution of Differential Problems. {\em Axioms} {\bf 7}(2) (2018) article n.\,36. \url{http://dx.doi.org//10.3390/axioms7020036}

\bibitem{BIT2009} L.\,Brugnano, F.\,Iavernaro, D.\,Trigiante. Hamiltonian BVMs (HBVMs): A family of ``drift-free'' methods for integrating polynomial Hamiltonian systems. {\em AIP Conf. Proc.} {\bf 1168} (2009) 715--718.

\bibitem{BIT2010} L.\,Brugnano, F.\,Iavernaro, D.\,Trigiante.  Hamiltonian Boundary Value Methods (Energy Preserving Discrete Line Integral Methods).  {\em JNAIAM J. Numer. Anal. Ind. Appl. Math.} {\bf 5},\,1-2 (2010) 17--37.

\bibitem{BIT2011} L.\,Brugnano, F.\,Iavernaro, D.\,Trigiante. A note on the efficient implementation of Hamiltonian BVMs. {\em J. Comput. Appl. Math.} {\bf 236} (2011) 375--383.

\bibitem{BIT2012} L.\,Brugnano, F.\,Iavernaro, D.\,Trigiante.  A simple framework for the derivation and analysis of effective one-step methods for ODEs. {\em Appl. Math. Comput.} {\bf 218} (2012) 8475--8485.

\bibitem{BIT2012_1} L.\,Brugnano, F.\,Iavernaro, D.\,Trigiante. A two-step, fourth-order method with energy preserving properties. {\em  Comput. Phys. Commun.} {\bf 183} (2012) 1860--1868.

\bibitem{BIT2012_2} L.\,Brugnano, F.\,Iavernaro, D.\,Trigiante. Energy and QUadratic Invariants Preserving integrators based upon Gauss collocation formulae. {\em SIAM J. Numer. Anal.} {\bf 50}, No.\,6 (2012) 2897--2916.

\bibitem{BIT2015}  L.\,Brugnano, F.\,Iavernaro, D.\,Trigiante. Analisys of Hamiltonian Boundary Value Methods (HBVMs): A class of energy-preserving Runge-Kutta methods for the numerical solution of polynomial Hamiltonian systems. {\em Commun. Nonlinear Sci. Numer. Simul.} {\bf 20} (2015) 650--667. 

\bibitem{BM2002} L.\,Brugnano, C.\,Magherini. Blended Implementation of Block Implicit Methods for ODEs. {\em Appl. Numer. Math.} {\bf 42} (2002) 29--45.

\bibitem{BM2007} L.\,Brugnano, C.\,Magherini. Blended Implicit Methods for solving ODE and DAE problems, and their extension for second order problems. {\em J. Comput. Appl. Math.} {\bf 205} (2007) 777--790.

\bibitem{EBQ2015} C.L.\,Ellison, J.W.\,Burby, H.\,Qin. Comment on ``Symplectic integration of magnetic systems'': A proof that the Boris algorithm is not variational. {\em J. Comput. Phys.} {\bf 301} (2015) 48--493.

\bibitem{G1999} D.J.\,Griffiths. {\em Introduction to electrodynamics (3rd ed.)}. Prentice-Hall, Upper Saddle River, NJ., 1999.

\bibitem{HL2018} E.\,Hairer, C.\,Lubich. Energy behaviour of the Boris method for charged-particle dynamics. {\em BIT} {\bf 58} (2018) 969--979.

\bibitem{HSLQ2015} Y.\,He, Y.\,Sun, J.\,Liu, H.\,Qin. Volume-preserving algorithms for charged particle dynamics. {\em J.\, Comput. Physics} {\bf 281} (2015) 135--147.

\bibitem{HSLQ2016} Y.\,He, Y.\,Sun, J.\,Liu, H.\,Qin. Higher order volume-preserving schemes for charged particle dynamics. {\em J.\, Comput. Physics} {\bf 305} (2016) 172--184.

\bibitem{HZSLQ2017} Y.\,He, Z.\,Zhou, Y.\,Sun, J.\,Liu, H.\,Qin. Explicit $K$-symplectic algorithms for charged particle dynamics. {\em Phys. Letters A} {\bf 381} (2017) 568--573.

\bibitem{IaPa2007} F.\,Iavernaro, B.\,Pace. $s$-stage trapezoidal methods for the conservation of Hamiltonian functions of polynomial type. {\em AIP Conf. Proc.} {\bf 936} (2007) 603--606.

\bibitem{IaPa2008} F.\,Iavernaro, B.\,Pace. Conservative block-Boundary Value Methods for the solution of polynomial Hamiltonian systems. {\em AIP Conf. Proc.} {\bf 1048} (2008) 888--891.

\bibitem{IaTr2009} F.\,Iavernaro, D.\,Trigiante. High-order Symmetric Schemes for the Energy Conservation of Polynomial Hamiltonian Problems.  {\em  JNAIAM. J. Numer. Anal. Ind. Appl. Math.} {\bf 4}, No.\,1-2 (2009) 87--101.

\bibitem{LW2016} H.\,Li, Y.\,Wang. A discrete line integral method of order two for the Lorentz force system. {\em Appl. Math. Comput.} {\bf 291} (2016) 207--212.

\bibitem{QZXLST2013} H.\,Qin, S.X.\,Zhang, J.Y.\,Xiao, J.\,Liu, Y.J.\,Sun, W.M.\,Tang. Why is Boris algorithm so good? {\em Physics of Plasmas} {\bf 20} (2013) 084503.

\bibitem{T2016} M.\,Tao. Explicit high-order symplectic integrators for charged particles in general electromagnetic fields. {\em J.\, Comput. Physics} {\bf 327} (2016) 245--251.
\bibitem{U2018} T.\,Umeda. A three-step Boris integrator for Lorentz force equation of charged particles. {\em  Comput. Phys. Commun.} {\bf 228} (2018) 1--4.

\bibitem{U2019} T.\,Umeda. Multi-step Boris rotation schemes for Lorentz force equation of charged particles. {\em  Comput. Phys. Commun.} {\bf 237} (2019) 37--41.

\end{thebibliography}
\end{document}